\input amstex
\documentstyle{amsppt}
%----------------------------------------------------------------
% Title:     Commutation relationships and curvature spin-tensors 
%            for extended spinor connections.
% Authors:   Ruslan Sharipov
% Comments:  AmSTeX, 22 pages, amsppt style
% MSC-class: 53A45, 53B15, 53C27, 83C60
%----------------------------------------------------------------
%           Replacement for output macro definition
%
\catcode`@=11
\redefine\output@{%
  \def\break{\penalty-\@M}\let\par\endgraf
  \ifodd\pageno\global\hoffset=105pt\else\global\hoffset=8pt\fi  
  \shipout\vbox{%
    \ifplain@
      \let\makeheadline\relax \let\makefootline\relax
    \else
      \iffirstpage@ \global\firstpage@false
        \let\rightheadline\frheadline
        \let\leftheadline\flheadline
      \else
        \ifrunheads@ %\let\makefootline\relax
        \else \let\makeheadline\relax
        \fi
      \fi
    \fi
    \makeheadline \pagebody \makefootline}%
  \advancepageno \ifnum\outputpenalty>-\@MM\else\dosupereject\fi
}
\catcode`\@=\active
%----------------------------------------------------------------
\nopagenumbers
\def\negskp{\hskip -2pt}
\def\Ker{\operatorname{Ker}}
\def\MatGrSL{\operatorname{SL}}
\def\MatGrSO{\operatorname{SO}}
\def\Alpha{\operatorname{A}}
\def\msum#1{\operatornamewithlimits{\sum^#1\!{\ssize\ldots}\!\sum^#1}}
\def\msums#1#2{\operatornamewithlimits{\sum^#1\!{\ssize\ldots}\!\sum^#2}}
\accentedsymbol\tx{\tilde x}
\accentedsymbol\vnabla{\nabla\kern -7pt\raise 5pt\vbox{\hrule width 5.5pt}
\kern 1.7pt}
\accentedsymbol\bPsi{\kern 1pt\overline{\kern -1pt\boldsymbol\Psi
\kern -1pt}\kern 1pt}
\def\blue#1{#1}
\catcode`#=11\def\diez{#}\catcode`#=6
\catcode`_=11\def\podcherkivanie{_}\catcode`_=8
\catcode`~=11\catcode`~=\active
\def\mycite#1{\cite{\blue{#1}}\immediate\special{ps:
     ShrHPSdict begin /ShrBORDERthickness 0 def}}

\def\mytag#1{%
    \tag#1}
\def\mythetag#1{\thetag{\blue{#1}}\immediate\special{ps:
     ShrHPSdict begin /ShrBORDERthickness 0 def}}
\def\myrefno#1{\no#1}
\def\myhref#1#2{\blue{#2}\immediate\special{ps:
     ShrHPSdict begin /ShrBORDERthickness 0 def}}
\def\myEarXivlink{\myhref{http://arXiv.org}{http:/\negskp/arXiv.org}}
\def\myGeoCities{\myhref{http://www.geocities.com}{GeoCities}}
\def\mytheorem#1{\csname proclaim\endcsname{Theorem #1}}
\def\mythetheorem#1{\blue{#1}\immediate\special{ps:
     ShrHPSdict begin /ShrBORDERthickness 0 def}}
\def\mylemma#1{\csname proclaim\endcsname{Lemma #1}}

\def\mydefinition#1{\definition{Definition #1}}
\def\mythedefinition#1{\blue{#1}\immediate\special{ps:
     ShrHPSdict begin /ShrBORDERthickness 0 def}}

%----------------------------------------------------------------
% Cyrillic fonts definition
%\font\tencyr=wncyr10
%\font\sevencyr=wncyr7
%\font\sevencyi=wncyi7
%----------------------------------------------------------------
\pagewidth{360pt}
\pageheight{606pt}
\topmatter
\title
Commutation relationships and curvature spin-tensors 
for extended spinor connections.
\endtitle
\author
R.~A.~Sharipov
\endauthor
\address 5 Rabochaya street, 450003 Ufa, Russia\newline
\vphantom{a}\kern 12pt Cell Phone: +7-(917)-476-93-48
\endaddress
\email \vtop to 30pt{\hsize=280pt\noindent
\myhref{mailto:R\podcherkivanie Sharipov\@ic.bashedu.ru}
{R\_\hskip 1pt Sharipov\@ic.bashedu.ru}\newline
\myhref{mailto:r-sharipov\@mail.ru}
{r-sharipov\@mail.ru}\newline
\myhref{mailto:ra\podcherkivanie sharipov\@lycos.com}{ra\_\hskip 1pt
sharipov\@lycos.com}\vss}
\endemail
\urladdr
\vtop to 20pt{\hsize=280pt\noindent
\myhref{http://www.geocities.com/r-sharipov}
{http:/\negskp/www.geocities.com/r-sharipov}\newline
\myhref{http://www.freetextbooks.boom.ru/index.html}
{http:/\negskp/www.freetextbooks.boom.ru/index.html}\vss}
\endurladdr
\abstract
    Extended spinor connections associated with composite 
spin-tensorial bundles are considered. Commutation relationships 
for covariant and multivariate differentiations and corresponding
curvature spin-tensors are derived. 
\endabstract
\endtopmatter
\loadbold
\loadeufb
\TagsOnRight
\document
%\input countstyle
%\vskip -20pt

%\special{header=resource.eps}
\rightheadtext{Commutation relationships and curvature spin-tensors \dots}
\head
1. Introduction of a spinor bundle.
\endhead
    This paper is a continuation of the paper \mycite{1}. For this reason
we keep all notations used in \mycite{1} and do not give any historical
background referring the reader backward to the paper \mycite{1} and to the
papers prior to it.\par
    Let $M$ be the {\it space-time} manifold, i\.\,e\. this is a 
$4$-dimensional orientable manifold equipped with a pseudo-Euclidean metric
$\bold g$ of the Minkowski-type signature $(+,-,-,-)$ and carrying a special
smooth geometric structure which is called a {\it polarization}. Once some
polarization is fixed, one can distinguish the {\it Future light cone\/}
from the {\it Past light cone\/} at each point $p\in M$ (see \mycite{2} for
more details). Moreover, we assume that $M$ admits the spinor structure.
This means that there is a two-dimensional smooth complex vector bundle
$SM$ over $M$ equipped with a skew-symmetric spin-tensorial field $\bold d$.
This spin-tensorial field $\bold d$ is called the {\it spin-metric tensor},
while $SM$ is called the {\it spinor bundle}.\par
     The spinor bundle $SM$ differs from a general two-dimensional complex
vector bundle over $M$ by its close relation to the tangent bundle $TM$.
The spin-metric tensor $\bold d$ and the metric tensor $\bold g$ are two
basic structures establishing this relation. Any local trivialization of
a two-dimensional vector bundle is given by two smooth sections of this
bundle $\boldsymbol\Psi_1$ and $\boldsymbol\Psi_2$ which are 
$\Bbb C$-linearly independent at each point $p$ of some open set $U\subset
M$. These two spinor fields $\boldsymbol\Psi_1$ and $\boldsymbol\Psi_2$ form
a moving frame $(U,\,\boldsymbol\Psi_1,\,\boldsymbol\Psi_2)$. A moving
frame $(U,\,\boldsymbol\Psi_1,\,\boldsymbol\Psi_2)$ is called an
orthonormal frame if
$$
\hskip -2em
d_{ij}=d(\boldsymbol\Psi_i,\boldsymbol\Psi_j)
=\Vmatrix 0 & 1\\ 
\vspace{1ex} -1 & 0\endVmatrix,
\mytag{1.1}
$$
i\.\,e\. if the spin metric tensor $\bold d$ is given by the skew-symmetric
matrix \mythetag{1.1} in this frame. Similarly, a moving frame
$(U,\,\boldsymbol\Upsilon_0,\,\boldsymbol\Upsilon_1,\,\boldsymbol\Upsilon_2,
\,\boldsymbol\Upsilon_3)$ of the tangent bundle $TM$ is given by four
smooth vector fields $\boldsymbol\Upsilon_0$, $\boldsymbol\Upsilon_1$, 
$\boldsymbol\Upsilon_2$, $\boldsymbol\Upsilon_3$ which are
$\Bbb R$-linearly independent at each point $p$ of the open set $U\subset
M$. \pagebreak This moving frame is called a positively polarized right
orthonormal frame if the following conditions are fulfilled:
\roster
\rosteritemwd=1pt
\item the value of the first vector filed $\boldsymbol\Upsilon_0$ at each
point $p\in U$ belongs to the interior of the Future light cone determined 
by the polarization of $M$;
\item it is a right frame in the sense of the orientation of $M$;
\item the metric tensor $\bold g$ is given by the standard
Minkowski matrix in this frame:
$$
\hskip -2em
g_{ij}=g(\boldsymbol\Upsilon_i,\boldsymbol\Upsilon_j)=\Vmatrix
\format \l&\quad\r&\quad\r&\quad\r\\
1 &0 &0 &0\\ 0 &-1 &0 &0\\ 0 &0 &-1 &0\\ 0 &0 &0 &-1\endVmatrix.
\mytag{1.2}
$$
\endroster\par
     If we have two orthonormal moving frames $(U,\,\boldsymbol\Psi_1,
\,\boldsymbol\Psi_2)$ and $(\tilde U,\,\tilde{\boldsymbol\Psi}_1,\,
\tilde{\boldsymbol\Psi}_2)$ of the spinor bundle $SM$ with overlapping 
domains $U\cap\tilde U\neq\varnothing$, then at each point 
$p\in U\cap\tilde U$ we have the following transition formulas:
$$
\xalignat 2
&\hskip -2em
\tilde{\boldsymbol\Psi}_i=\sum^2_{j=1}\goth S^j_i\,
\boldsymbol\Psi_j,
&&\boldsymbol\Psi_i=\sum^2_{j=1}\goth T^j_i\,
\tilde{\boldsymbol\Psi}_j,
\mytag{1.3}
\endxalignat
$$
The transition matrices $\goth S$ and $\goth T$ in \mythetag{1.3} are
inverse to each other: $\goth T=\goth S^{-1}$. From \mythetag{1.1} for
these matrices one easily derives
$$
\xalignat 2
&\hskip -2em
\goth S(p)\in\MatGrSL(2,\Bbb C),
&&\goth T(p)\in\MatGrSL(2,\Bbb C).
\mytag{1.4}
\endxalignat
$$
In a similar way, if we have two positively polarized right orthonormal
frames $(U,\,\boldsymbol\Upsilon_0,\,\boldsymbol\Upsilon_1,\,\boldsymbol
\Upsilon_2,\,\boldsymbol\Upsilon_3)$ and $(\tilde U,\,
\tilde{\boldsymbol\Upsilon}_0,\,\tilde{\boldsymbol\Upsilon}_1,\,
\tilde{\boldsymbol\Upsilon}_2,\,\tilde{\boldsymbol\Upsilon}_3)$ of the
tangent bundle $TM$ with overlapping domains $U\cap\tilde U\neq\varnothing$,
then we have
$$
\xalignat 2
&\hskip -2em
\tilde{\boldsymbol\Upsilon}_i=\sum^3_{j=0}S^j_i\,\boldsymbol\Upsilon_j,
&&\boldsymbol\Upsilon_i=\sum^3_{j=0}T^j_i\,\tilde{\boldsymbol\Upsilon}_j.
\mytag{1.5}
\endxalignat
$$
The matrices $S$ and $T$ in the formulas \mythetag{1.5} are inverse to each
other: $T=S^{-1}$. From \mythetag{1.2} and from the above conditions
\therosteritem{1}--\therosteritem{3} for each point $p\in U\cap\tilde U$ we
get
$$
\xalignat 2
&S(p)\in\MatGrSO^+(1,3,\Bbb R),
&&T(p)\in\MatGrSO^+(1,3,\Bbb R).
\endxalignat
$$
Note that the special linear group $\MatGrSL(2,\Bbb C)$ in \mythetag{1.4}
and the special orthochronous Lorentzian group $\MatGrSO^+(1,3,\Bbb R)$
are related by the canonical homomorphism 
$$
\hskip -2em
\varphi\!:\,\MatGrSL(2,\Bbb C)\to\MatGrSO^+(1,3,\Bbb R).
\mytag{1.6}
$$
The canonical homomorphism $\varphi$ in \mythetag{1.6} is a surjective
mapping. Its kernel is a discrete group. It is composed by the following 
two matrices:
$$
\xalignat 2
&\hskip -2em
\boldsymbol\sigma_0=\Vmatrix 1 & 0\\0 & 1\endVmatrix,
&&-\boldsymbol\sigma_0=\Vmatrix -1 & 0\\0 & -1\endVmatrix.
\mytag{1.7}
\endxalignat
$$
Since $\Ker\varphi$ is a discrete set composed by two matrices 
\mythetag{1.7}, topologically $\varphi$ is a double sheeted not
ramified covering of real analytic manifolds. It can be given
by explicit formulas in terms of the matrix components. If
$$
S=\varphi(\goth S),
$$
then for the components of the matrix $S$ we have the following 
expressions:
$$
\gather
\hskip -6em
\gathered
S^0_0=\frac{\overline{\goth S^1_1}\,\goth S^1_1
+\overline{\goth S^1_2}\,\goth S^1_2+\overline{\goth S^2_1}
\,\goth S^2_1+\overline{\goth S^2_2}\,\goth S^2_2}{2},\\
S^0_1=\frac{\overline{\goth S^1_1}\,\goth S^1_2
+\overline{\goth S^1_2}\,\goth S^1_1+\overline{\goth S^2_1}
\,\goth S^2_2+\overline{\goth S^2_2}\,\goth S^2_1}{2},\\
S^0_2=\frac{\overline{\goth S^1_2}\,\goth S^1_1
-\overline{\goth S^1_1}\,\goth S^1_2+\overline{\goth S^2_2}
\,\goth S^2_1-\overline{\goth S^2_1}\,\goth S^2_2}{2\,i},\\
S^0_3=\frac{\overline{\goth S^1_1}\,\goth S^1_1
-\overline{\goth S^1_2}\,\goth S^1_2+\overline{\goth S^2_1}
\,\goth S^2_1-\overline{\goth S^2_2}\,\goth S^2_2}{2},
\endgathered
\mytag{1.8}\\
\vspace{1ex}
\hskip 2em
\gathered
S^1_0=\frac{\overline{\goth S^2_1}\,\goth S^1_1
+\overline{\goth S^1_1}\,\goth S^2_1+\overline{\goth S^2_2}
\,\goth S^1_2+\overline{\goth S^1_2}\,\goth S^2_2}{2},\\
S^1_1=\frac{\overline{\goth S^1_2}\,\goth S^1_2
+\overline{\goth S^2_1}\,\goth S^2_1+\overline{\goth S^2_2}
\,\goth S^1_1+\overline{\goth S^1_1}\,\goth S^2_2}{2},\\
S^1_2=\frac{\overline{\goth S^1_2}\,\goth S^2_1
-\overline{\goth S^2_1}\,\goth S^1_2+\overline{\goth S^2_2}
\,\goth S^1_1-\overline{\goth S^1_1}\,\goth S^2_2}{2\,i},\\
S^1_3=\frac{\overline{\goth S^2_1}\,\goth S^1_1
+\overline{\goth S^1_1}\,\goth S^2_1-\overline{\goth S^2_2}
\,\goth S^1_2-\overline{\goth S^1_2}\,\goth S^2_2}{2},
\endgathered
\mytag{1.9}\\
\vspace{1ex}
\hskip -6em
\gathered
S^2_0=\frac{\overline{\goth S^1_1}\,\goth S^2_1
-\overline{\goth S^2_1}\,\goth S^1_1+\overline{\goth S^1_2}
\,\goth S^2_2-\overline{\goth S^2_2}\,\goth S^1_2}{2\,i},\\
S^2_1=\frac{\overline{\goth S^1_2}\,\goth S^2_1
-\overline{\goth S^2_1}\,\goth S^1_2+\overline{\goth S^1_1}
\,\goth S^2_2-\overline{\goth S^2_2}\,\goth S^1_1}{2\,i},\\
S^2_2=\frac{\overline{\goth S^2_2}\,\goth S^1_1
+\overline{\goth S^1_1}\,\goth S^2_2-\overline{\goth S^2_1}
\,\goth S^1_2-\overline{\goth S^1_2}\,\goth S^2_1}{2},\\
S^2_3=\frac{\overline{\goth S^1_1}\,\goth S^2_1
-\overline{\goth S^2_1}\,\goth S^1_1+\overline{\goth S^2_2}
\,\goth S^1_2-\overline{\goth S^1_2}\,\goth S^2_2}{2\,i},
\endgathered
\mytag{1.10}\\
\vspace{1ex}
\hskip 2em
\gathered
S^3_0=\frac{\overline{\goth S^1_1}\,\goth S^1_1
+\overline{\goth S^1_2}\,\goth S^1_2-\overline{\goth S^2_1}
\,\goth S^2_1-\overline{\goth S^2_2}\,\goth S^2_2}{2},\\
S^3_1=\frac{\overline{\goth S^1_1}\,\goth S^1_2
+\overline{\goth S^1_2}\,\goth S^1_1-\overline{\goth S^2_1}
\,\goth S^2_2-\overline{\goth S^2_2}\,\goth S^2_1}{2},\\
S^3_2=\frac{\overline{\goth S^1_2}\,\goth S^1_1
-\overline{\goth S^1_1}\,\goth S^1_2+\overline{\goth S^2_1}
\,\goth S^2_2-\overline{\goth S^2_2}\,\goth S^2_1}{2\,i},\\
S^3_3=\frac{\overline{\goth S^1_1}\,\goth S^1_1
+\overline{\goth S^2_2}\,\goth S^2_2-\overline{\goth S^2_1}
\,\goth S^2_1-\overline{\goth S^1_2}\,\goth S^1_2}{2}.
\endgathered
\mytag{1.11}
\endgather
$$
The first formula \mythetag{1.8} was presented in \mycite{1}. The
whole set of the above formulas \mythetag{1.8}, \mythetag{1.9},
\mythetag{1.10}, and \mythetag{1.11} can be found in \mycite{3}.
\mydefinition{1.1} Let $SM$ be a two-dimensional complex vector 
bundle over the space-time manifold $M$ equipped with a nonzero
spin-metric $\bold d$. It is called the {\it spinor bundle} if each
orthonormal frame $(U,\,\boldsymbol\Psi_1,\,\boldsymbol\Psi_2)$ of
$SM$ is associated with some positively polarized right orthonormal 
frame $(U,\,\boldsymbol\Upsilon_0,\,\boldsymbol\Upsilon_1,\,
\boldsymbol\Upsilon_2,\,\boldsymbol\Upsilon_3)$ of the tangent bundle
$TM$ such that for any two orthonormal frames $(U,\,\boldsymbol\Psi_1,
\,\boldsymbol\Psi_2)$ and $(\tilde U,\,\tilde{\boldsymbol\Psi}_1,\,
\tilde{\boldsymbol\Psi}_2)$ with overlapping domains $U\cap\tilde U\neq
\varnothing$ the associated tangent frames $(U,\,\boldsymbol\Upsilon_0,
\,\boldsymbol\Upsilon_1,\,\boldsymbol\Upsilon_2,\,\boldsymbol\Upsilon_3)$
and $(\tilde U,\,\tilde{\boldsymbol\Upsilon}_0,\,
\tilde{\boldsymbol\Upsilon}_1,\,\tilde{\boldsymbol\Upsilon}_2,\,
\tilde{\boldsymbol\Upsilon}_3)$ are related to each other by means of
the formulas \mythetag{1.5}, where the transition matrices $S$ and $T$
are obtained from the transition matrices $\goth S$ and $\goth T$ in
\mythetag{1.3} by applying the group homomorphism \mythetag{1.6}, 
i\.\,e\. $S=\varphi(\goth S)$ and $T=\varphi(\goth T)$.
\enddefinition
    Let $(U,\,\boldsymbol\Psi_1,\,\boldsymbol\Psi_2)$ be an orthonormal
frame of the spinor bundle $SM$ and assume that the domain $U$ is
sufficiently small to be equipped with some local coordinates $x^0,\,
x^1,\,x^2,\,x^3$. Then $U$ is a local chart and we have the holonomic
frame in $TM$ determined by the local coordinates of this chart:
$$
\hskip -2em
\bold E_0=\frac{\partial}{\partial x^0},\quad 
\bold E_1=\frac{\partial}{\partial x^1},\quad 
\bold E_2=\frac{\partial}{\partial x^2},\quad 
\bold E_3=\frac{\partial}{\partial x^3}.
\mytag{1.12}
$$
Passing from $(U,\,\boldsymbol\Psi_1,\,\boldsymbol\Psi_2)$ to the
associated frame $(U,\,\boldsymbol\Upsilon_0,\,\boldsymbol\Upsilon_1,\,
\boldsymbol\Upsilon_2,\,\boldsymbol\Upsilon_3)$, in general case we find
that it doesn't coincide with the coordinate frame $(U,\,\bold E_0,
\,\bold E_1,\,\bold E_2,\,\bold E_3)$ since in general case 
$(U,\,\boldsymbol\Upsilon_0,\,\boldsymbol\Upsilon_1,\,
\boldsymbol\Upsilon_2,\,\boldsymbol\Upsilon_3)$ is a non-holonomic frame.
Let's consider the following expansion of the non-holonomic frame vectors:
$$
\hskip -2em
\boldsymbol\Upsilon_i=\sum^3_{j=0}\Upsilon^j_i\,\bold E_j.
\mytag{1.13}
$$
In the case of a holonomic frame \mythetag{1.12} all mutual commutators
of the vector fields $\bold E_0,\,\bold E_1,\,\bold E_2,\,\bold E_3$ are
equal to zero. In the case of a non-holonomic frame it is not so:
$$
\hskip -2em
[\boldsymbol\Upsilon_i,\,\boldsymbol\Upsilon_j]=
\sum^3_{k=0}c^{\,k}_{ij}\,\boldsymbol\Upsilon_k.
\mytag{1.14}
$$
From \mythetag{1.13} and  \mythetag{1.14} one easily derives
$$
\hskip -2em
\sum^3_{s=0}\Upsilon^s_i\,\frac{\partial\Upsilon^m_j}{\partial x^s}
-\sum^3_{s=0}\Upsilon^s_j\,\frac{\partial\Upsilon^m_i}{\partial x^s}
=\sum^3_{k=0}c^{\,k}_{ij}\,\Upsilon^m_k.
\mytag{1.15}
$$
The coefficients $c^{\,k}_{ij}$ can be calculated using either 
\mythetag{1.14} or \mythetag{1.15}. They form a frame specific set
of functions.
\head
2. Tensors and spin-tensors.
\endhead
    Let $M$ be the space-time and let $p$ be a point of $M$. Then
$T_p(M)$ is a tangent space of $M$ at the point $p$. Similarly, 
$T^*_p(M)$ is a cotangent space at the same point, it is dual to
the space $T_p(M)$. Both $T_p(M)$ and $T^*_p(M)$ are real linear
vector spaces. Following the recipe of \mycite{1}, we introduce
their complexifications:
$$
\xalignat 2
&\hskip -2em
\Bbb CT_p(M)=\Bbb C\otimes T_p(M),
&&\Bbb CT^*_p(M)=\Bbb C\otimes T^*_p(M).
\mytag{2.1}
\endxalignat
$$
Complexified tensor spaces then are introduced as multiple tensor 
products of several copies of the spaces $\Bbb CT_p(M)$ and 
$\Bbb CT^*_p(M)$ introduced in \mythetag{2.1}:
$$
\Bbb CT^m_n(p,M)=\overbrace{\Bbb CT_p(M)\otimes\ldots\otimes
\Bbb CT_p(M)}^{\text{$m$ times}}\otimes\underbrace{\Bbb CT^*_p(M)
\otimes\ldots\otimes \Bbb CT^*_p(M)}_{\text{$n$ times}}.
\quad
\mytag{2.2}
$$
The spinor bundle $SM$ is a complex vector bundle over $M$ from the
very beginning. Let's denote by $S_p(M)$ its fiber over the point 
$p\in M$. Let $S^*_p(M)$ be its dual space. Moreover, we consider
the hermitian conjugate space $S^{\sssize\dagger}_p(M)$ for $S_p(M)$
and its dual space $S^{{\sssize\dagger}*}_p(M)=S^{*\sssize\dagger}_p(M)$.
Then we can define the following tensor products:
$$
\align
&\hskip -2em
S^\alpha_\beta(p,M)=\overbrace{S_p(M)\otimes\ldots\otimes 
S_p(M)}^{\text{$\alpha$ times}}\otimes
\underbrace{S^*_p(M)\otimes\ldots\otimes 
S^*_p(M)}_{\text{$\beta$ times}},
\mytag{2.3}\\
&\hskip -2em
\bar S^\nu_\gamma(p,M)=\overbrace{S^{{\sssize\dagger}*}_p(M)\otimes
\ldots\otimes S^{{\sssize\dagger}*}_p(M)}^{\text{$\nu$ times}}
\otimes\underbrace{S^{\sssize\dagger}_p(M)\otimes\ldots\otimes 
S^{\sssize\dagger}_p(M)}_{\text{$\gamma$ times}}.
\mytag{2.4}
\endalign
$$     
Combining \mythetag{2.2}, \mythetag{2.3}, and \mythetag{2.4}, we define
one more tensor product
$$
\hskip -2em
S^\alpha_\beta\bar S^\nu_\gamma T^r_s(p,M)=S^\alpha_\beta(p,M)
\otimes\bar S^\nu_\gamma(p,M)\otimes\Bbb CT^m_n(p,M).
\mytag{2.5}
$$
Elements of the space \mythetag{2.5} are called {\it spin-tensors} of the
type $(\alpha,\beta|\nu,\gamma|m,n)$ at the point $p\in M$. Elements of 
other three spaces \mythetag{2.2}, \mythetag{2.3}, and \mythetag{2.4}
are also called spin-tensors, though they are special cases of a
general spin-tensorial object. The spin-tensorial space \mythetag{2.5}
admits the canonical semilinear isomorphism $\tau$:
$$
\hskip -2em
\tau\!:\,
S^\alpha_\beta\bar S^\nu_\gamma T^m_n(p,M)
\to S^\nu_\gamma\bar S^\alpha_\beta T^m_n(p,M)
\mytag{2.6}
$$
(see the definition of $\tau$ and more details concerning it in
\mycite{1}). The spin-tensorial spaces \mythetag{2.5} with $p$ running
over $M$ are glued into a bundle. It is called the {\it spin-tensorial
bundle} of the type $(\alpha,\beta|\nu,\gamma|m,n)$ and is denoted by
$S^\alpha_\beta\bar S^\nu_\gamma T^m_n\!M$. Then the isomorphisms
\mythetag{2.6} with the point $p$ running over $M$ are glued into a
semilinear isomorphism of spin-tensorial bundles:
$$
\tau\!:\,S^\alpha_\beta\bar S^\nu_\gamma T^m_n\!M
\to S^\nu_\gamma\bar S^\alpha_\beta T^m_n\!M.
$$
A traditional spin-tensorial field of the type $(\alpha,\beta|\nu,
\gamma|m,n)$ by definition  is a local or global smooth section of
the spin-tensorial bundle $S^\alpha_\beta\bar S^\nu_\gamma T^m_n\!M$.
Non-traditional (extended) spin-tensorial fields were introduced in
\mycite{1} along with non-traditional (extended) connections. We 
shall give their definitions a little bit later.
\head
3. Coordinate representation of spin-tensorial fields.
\endhead
     In order to represent a tensorial field in a coordinate form
it is sufficient to have a local chart in $M$. In the case of spin-tensorial
fields, in addition to a local chart, we need to have two frames: one in
$SM$ and the other in $TM$. 
\mydefinition{3.1} Let $U$ be a local chart of the space-time manifold
$M$. We say that $U$ is an {\it equipped local chart} if there is an
orthonormal spinor frame $(U,\,\boldsymbol\Psi_1,\,\boldsymbol\Psi_2)$
with the domain $U$ and, hence, according to the 
definition~\mythedefinition{1.1}, there is a positively polarized right
orthonormal tangent frame $(U,\,\boldsymbol\Upsilon_0,\,
\boldsymbol\Upsilon_1,\,\boldsymbol\Upsilon_2,\,\boldsymbol\Upsilon_3)$
canonically associated with the frame $(U,\,\boldsymbol\Psi_1,
\,\boldsymbol\Psi_2)$.
\enddefinition
     Equipped local charts cover the whole space-time manifold, i\.\,e\.
they form an atlas. Therefore, they describe completely the structure 
of the the space-time manifold $M$ and its bundles $SM$ and $TM$.\par
     Let $U$ be an equipped local chart with the local coordinates 
$x^0,\,x^1,\,x^2,\,x^3$. Assume that $(U,\,\boldsymbol\Psi_1,\,
\boldsymbol\Psi_2)$ and $(U,\,\boldsymbol\Upsilon_0,\,
\boldsymbol\Upsilon_1,\,\boldsymbol\Upsilon_2,\,\boldsymbol\Upsilon_3)$
are two frames associated with $U$ and thus being its equipment. Denote
by $(U,\,\boldsymbol\vartheta^{\,1},\,\boldsymbol\vartheta^{\,2})$ the 
dual cospinor frame for $(U,\,\boldsymbol\Psi_1,\,\boldsymbol\Psi_2)$ and
denote by $(U,\,\boldsymbol\eta^0,\,\boldsymbol\eta^1,\,\boldsymbol\eta^2,\,
\boldsymbol\eta^3)$ the dual covectorial frame for $(U,\,
\boldsymbol\Upsilon_0,\,\boldsymbol\Upsilon_1,\,\boldsymbol\Upsilon_2,\,
\boldsymbol\Upsilon_3)$. Moreover, we denote
$$
\xalignat 2
&\hskip -2em
\bPsi_i=\tau(\boldsymbol\Psi_i),
&&\overline{\boldsymbol\vartheta}\vphantom{\boldsymbol\vartheta}^{\,i}
=\tau(\boldsymbol\vartheta^{\,i}).
\mytag{3.1}
\endxalignat
$$
The barred spinor fields \mythetag{3.1} form two frames
$(U,\,\bPsi_1,\,\bPsi_2)$ and $(U,\,\overline{\boldsymbol\vartheta}
\vphantom{\boldsymbol\vartheta}^{\,1},\,\overline{\boldsymbol\vartheta}
\vphantom{\boldsymbol\vartheta}^{\,2})$ dual to each other. Using all
the above frames, we define the following fields:
$$
\align
&\hskip -2em
\boldsymbol\Upsilon^{k_1\ldots\,k_n}_{h_1\ldots\,h_m}
=\boldsymbol\Upsilon_{h_1}\otimes\ldots\otimes\boldsymbol\Upsilon_{h_m}
\otimes\boldsymbol\eta^{k_1}\otimes\ldots\otimes\boldsymbol\eta^{k_n},
\mytag{3.2}\\
\vspace{1ex}
&\hskip -2em
\boldsymbol\Psi^{j_1\ldots\,j_\beta}_{i_1\ldots\,i_\alpha}
=\boldsymbol\Psi_{i_1}\otimes\ldots\otimes\boldsymbol\Psi_{i_\alpha}
\otimes\boldsymbol\vartheta^{\,j_1}\otimes\ldots\otimes
\boldsymbol\vartheta^{\,j_\beta},
\mytag{3.3}\\
\vspace{1ex}
&\hskip -2em
\bPsi^{\bar j_1\ldots\,\bar j_\gamma}_{\bar i_1\ldots\,\bar i_\nu}
=\bPsi_{\bar i_1}\otimes\ldots\otimes\bPsi_{\bar i_\nu}\otimes
\overline{\boldsymbol\vartheta}\vphantom{\boldsymbol\vartheta}^{\,j_1}
\otimes\ldots\otimes\overline{\boldsymbol\vartheta}
\vphantom{\boldsymbol\vartheta}^{\,j_\gamma}.
\mytag{3.4}
\endalign
$$
Then, using \mythetag{3.2}, \mythetag{3.3}, and \mythetag{3.4}, we introduce
the following tensor product:
$$
\hskip -2em
\boldsymbol\Psi^{j_1\ldots\,j_\beta\,\bar j_1\ldots\,\bar j_\gamma
\,k_1\ldots\, k_n}_{i_1\ldots\,i_\alpha\,\bar i_1\ldots\,\bar i_\nu
\,h_1\ldots\,h_m}
=\boldsymbol\Psi^{j_1\ldots\,j_\beta}_{i_1\ldots\,i_\alpha}
\otimes\bPsi^{\bar j_1\ldots\,\bar j_\gamma}_{\bar i_1\ldots\,
\bar i_\nu}\otimes\boldsymbol\Upsilon^{k_1\ldots\,k_n}_{h_1\ldots\,h_m}.
\mytag{3.5}
$$
It is easy to see that the formula \mythetag{3.5} defines a series of
local spin-tensorial fields of the type $(\alpha,\beta|\nu,\gamma|m,n)$
with the domain $U$. Assume that $\bold X$ is an arbitrary spin-tensorial
field of the same type $(\alpha,\beta|\nu,\gamma|m,n)$. For this 
spin-tensorial field within the domain $U$ we can write the following 
expansion:
$$
\bold X=\msum{2}\Sb i_1,\,\ldots,\,i_\alpha\\ j_1,\,\ldots,\,j_\beta
\endSb\msum{2}\Sb\bar i_1,\,\ldots,\,\bar i_\nu\\ \bar j_1,\,\ldots,\,
\bar j_\gamma\endSb\msum{3}\Sb h_1,\,\ldots,\,h_m\\ k_1,\,\ldots,\,k_n
\endSb
X^{i_1\ldots\,i_\alpha\,\bar i_1\ldots\,\bar i_\nu\,h_1\ldots\,h_m}_{j_1
\ldots\,j_\beta\,\bar j_1\ldots\,\bar j_\gamma\,k_1\ldots\, k_n}\ 
\boldsymbol\Psi^{j_1\ldots\,j_\beta\,\bar j_1\ldots\,\bar j_\gamma\,
k_1\ldots\, k_n}_{i_1\ldots\,i_\alpha\,\bar i_1\ldots\,\bar i_\nu\,
h_1\ldots\,h_m}.\quad
\mytag{3.6}
$$
The coefficients $X^{i_1\ldots\,i_\alpha\,\bar i_1\ldots\,\bar i_\nu\,
h_1\ldots\,h_m}_{j_1\ldots\,j_\beta\,\bar j_1\ldots\,\bar j_\gamma\,k_1
\ldots\, k_n}$ in the expansion \mythetag{3.6} are functions of the local
coordinates $x^0,\,x^1,\,x^2,\,x^3$ of a point $p\in U$:
$$
\hskip -2em
X^{i_1\ldots\,i_\alpha\,\bar i_1\ldots\,\bar i_\nu\,h_1\ldots\,h_m}_{j_1
\ldots\,j_\beta\,\bar j_1\ldots\,\bar j_\gamma\,k_1\ldots\, k_n}=
X^{i_1\ldots\,i_\alpha\,\bar i_1\ldots\,\bar i_\nu\,h_1\ldots\,h_m}_{j_1
\ldots\,j_\beta\,\bar j_1\ldots\,\bar j_\gamma\,k_1\ldots\, k_n}(x^0,
x^1,x^2,x^3).
\mytag{3.7}
$$
These functions \mythetag{3.7} are called the {\it components} of a 
spin-tensorial field $X$ in an equipped local chart $U$, while the
expansion \mythetag{3.6} itself is the {\it coordinate representation}
of the field $\bold X$. 
\head
4. Composite spin-tensorial bundles
and extended\\spin-tensorial fields.
\endhead
     Let $\bold S[1],\,\ldots,\,\bold S[J]$ be some spin-tensorial 
fields of various types, e\.\,g\. we can denote by $(\alpha_P,
\beta_P|\nu_P,\gamma_P|m_P,n_P)$ the type of the $P$-th field
$\bold S[P]$ in the series $\bold S[1],\,\ldots,\,\bold S[J]$. In 
some cases one need to treat $\bold S[1],\,\ldots,\,\bold S[J]$
not as actual fields, but as independent variables. For example,
if we consider a physical field theory with the fields $\bold S[1],
\,\ldots,\,\bold S[J]$, then the Lagrange function 
$$
\hskip -2em
\Cal L=\Cal L(p,\bold S[1],\ldots,\bold S[J])
\mytag{4.1}
$$
of this field theory is a function of several spin-tensorial arguments
$\bold S[1],\,\ldots,\,\bold S[J]$ and of one point argument $p\in M$.
Composite spin-tensorial bundles were introduced in \mycite{1} for to
formalize the argument set of the function \mythetag{4.1}. In the present
case the composite tensor bundle $N$ is defined as the following direct 
sum:
$$
\hskip -2em
N=S^{\alpha_1}_{\beta_1}\!\bar S^{\nu_1}_{\gamma_1}
T^{m_1}_{n_1}\!M\oplus\ldots\oplus
S^{\alpha_J}_{\beta_J}\!\bar S^{\nu_J}_{\gamma_J}
T^{m_J}_{n_J}\!M.\quad
\mytag{4.2}
$$
By definition a point $q$ of the composite spin-tensorial bundle
\mythetag{4.2} is a list
$$
\hskip -2em
q=(p,\,\bold S[1],\,\ldots,\,\bold S[J]),
\mytag{4.3}
$$
where $p$ is a point of the space-time $M$, while $\bold S[1],\,\ldots,
\,\bold S[J]$ are spin-tensors of the types $(\alpha_1,\beta_1|\nu_1,
\gamma_1|m_1,n_1),\,\ldots,\,(\alpha_J,\beta_J|\nu_J,\gamma_J|m_J,n_J)$
at the point $p$. 
\mydefinition{4.1} Let $N$ be a composite spin-tensorial bundle over 
the space-time manifold $M$ in the sense of the formula \mythetag{4.2}.
An extended spin-tensorial field $\bold X$ of the type 
$(\varepsilon,\eta|\sigma,\zeta|e,f)$ is a spin-tensor-valued function 
in $N$ such that it takes each point $q\in N$ to some spin-tensor 
$\bold X(q)\in S^{\,\varepsilon}_\eta\bar S^\sigma_\zeta T^e_f(p,M)$, 
where $p=\pi(q)$ is the projection of the point $q$.
\enddefinition
     Let $U$ be an equipped local chart of the space-time manifold $M$
(see the definition~\mythedefinition{3.1} above). Then a point $p\in U$
is given by its coordinates 
$$
\hskip -2em
x^0,\,x^1,\,x^2,\,x^3,
\mytag{4.4}
$$
while spin-tensors $\bold S[1],\,\ldots,\,\bold S[J]$ at the point $p$
are given by their components referred to the frames $(U,\,
\boldsymbol\Psi_1,\,\boldsymbol\Psi_2)$ and $(U,\,\boldsymbol\Upsilon_0,
\,\boldsymbol\Upsilon_1,\,\boldsymbol\Upsilon_2,\,\boldsymbol\Upsilon_3)$:
$$
\hskip -2em
\gathered
S^{1\kern 1pt\ldots\,1\,1\kern 1pt\ldots\,1\,0\kern 1pt\ldots\,
0}_{1\kern 1pt\ldots\,1\,1\kern 1pt\ldots\,1\,0\kern 1pt\ldots\,0}[1],
\,\ldots,\,
S^{2\kern 1pt\ldots\,2\,2\kern 1pt\ldots\,2\,3\kern 1pt\ldots\,
3}_{2\kern 1pt\ldots\,2\,2\kern 1pt\ldots\,2\,3\kern 1pt\ldots\,3}[1],
\,\ldots\\
\kern 10em
\ldots,\,
S^{1\kern 1pt\ldots\,1\,1\kern 1pt\ldots\,1\,0\kern 1pt\ldots\,
0}_{1\kern 1pt\ldots\,1\,1\kern 1pt\ldots\,1\,0\kern 1pt\ldots\,0}[J],
\,\ldots,\,
S^{2\kern 1pt\ldots\,2\,2\kern 1pt\ldots\,2\,3\kern 1pt\ldots\,
3}_{2\kern 1pt\ldots\,2\,2\kern 1pt\ldots\,2\,3\kern 1pt\ldots\,3}[J].
\endgathered
\mytag{4.5}
$$
The quantities \mythetag{4.4} and \mythetag{4.5} form a complete set
of variables that can be used as local coordinates for a point $q$
in \mythetag{4.3}. For an extended spin-tensorial field $\bold X$
of the type $(\varepsilon,\eta|\sigma,\zeta|e,f)$ one can write the
expansion
$$
\bold X=\msum{2}\Sb i_1,\,\ldots,\,i_\varepsilon\\ j_1,\,\ldots,\,j_\eta
\endSb\msum{2}\Sb\bar i_1,\,\ldots,\,\bar i_\sigma\\ \bar j_1,\,\ldots,\,
\bar j_\zeta\endSb\msum{3}\Sb h_1,\,\ldots,\,h_e\\ k_1,\,\ldots,\,k_f
\endSb
X^{i_1\ldots\,i_\varepsilon\,\bar i_1\ldots\,\bar i_\sigma\,h_1\ldots
\,h_e}_{j_1\ldots\,j_\eta\,\bar j_1\ldots\,\bar j_\zeta\,k_1\ldots\,k_f}
\ \boldsymbol\Psi^{j_1\ldots\,j_\eta\,\bar j_1\ldots\,\bar j_\zeta\,k_1
\ldots\, k_f}_{i_1\ldots\,i_\varepsilon\,\bar i_1\ldots\,\bar i_\sigma\,
h_1\ldots\,h_e}\quad
\mytag{4.6}
$$
similar to \mythetag{3.6}. However, unlike to \mythetag{3.7}, the 
coefficients $X^{i_1\ldots\,i_\varepsilon\,\bar i_1\ldots\,\bar i_\sigma
\,h_1\ldots\,h_e}_{j_1\ldots\,j_\eta\,\bar j_1\ldots\,\bar j_\zeta\,k_1
\ldots\,k_f}$ in the expansion \mythetag{4.6} are functions of the whole
set of variables \mythetag{4.4} and \mythetag{4.5}. Thus, taking an
equipped local chart $U$ of $M$, we get a coordinate representation 
for points of the composite spin-tensorial bundle \mythetag{4.2} and for
extended spin-tensorial fields associated with it.\par
     Under a change of equipped local charts the coordinates \mythetag{4.4}
are transformed traditionally by means of the transition functions
$$
\xalignat 2
&\cases
\tx^0=\tx^1(x^0,\,\ldots,x^3),\\
. \ . \ . \ .\ . \ . \ . \ . \ . \ . \ 
. \ . \ . \ . \ .\\
\tx^3=\tx^n(x^0,\,\ldots,x^3).
\endcases
&&\cases
x^0=\tx^1(\tx^0,\,\ldots,\tx^3),\\
. \ . \ . \ .\ . \ . \ . \ . \ . \ . \ 
. \ . \ . \ . \ .\\
x^3=\tx^n(\tx^0,\,\ldots,\tx^3),
\endcases
\endxalignat
$$
while the coordinates \mythetag{4.5} are transformed as the components
of spin-tensors
$$
\align
&\hskip -4em
\left\{
\aligned
&\tilde S^{i_1\ldots\,i_\alpha\bar i_1\ldots\,\bar i_\nu
h_1\ldots\,h_m}_{j_1\ldots\,j_\beta\bar j_1\ldots\,\bar j_\gamma
k_1\ldots\,k_n}[P]
=\dsize\msum{2}\Sb a_1,\,\ldots,\,a_\alpha\\ b_1,\,\ldots,\,b_\beta\endSb
\dsize\msum{2}\Sb \bar a_1,\,\ldots,\,\bar a_\nu\\ 
\bar b_1,\,\ldots,\,\bar b_\gamma\endSb
\dsize\msum{3}
\Sb c_1,\,\ldots,\,c_m\\ d_1,\,\ldots,\,d_n\endSb
\goth T^{\,i_1}_{a_1}\ldots\,\goth T^{\,i_\alpha}_{a_\alpha}\,\times\\
&\kern 40pt
\times\,\goth S^{b_1}_{j_1}\ldots\,\goth S^{b_\beta}_{j_\beta}
\ \overline{\goth T^{\,\bar i_1}_{\bar a_1}}
\ldots\,\overline{\goth T^{\,\bar i_\nu}_{\bar a_\nu}}\ 
\ \overline{\goth S^{\,\bar b_1}_{\bar j_1}}
\ldots\,\overline{\goth S^{\,\bar b_\gamma}_{\bar j_\gamma}}\ 
T^{h_1}_{c_1}\ldots\,T^{h_m}_{c_m}\,\times\\
\vspace{1.5ex}
&\kern 60pt
\times\,S^{\,d_1}_{k_1}\ldots\,S^{\,d_n}_{k_n}\ 
S^{\,a_1\ldots\,a_\alpha\,\bar a_1\ldots\,\bar a_\nu\,
c_1\ldots\,c_m}_{\,b_1\ldots\,b_\beta\,\bar b_1\ldots\,\bar b_\gamma
\,d_1\ldots\,d_n}[P],
\endaligned
\right.
\mytag{4.7}\\
&\hskip -4em
\left\{
\aligned
&S^{i_1\ldots\,i_\alpha\bar i_1\ldots\,\bar i_\nu
h_1\ldots\,h_m}_{j_1\ldots\,j_\beta\bar j_1\ldots\,\bar j_\gamma
k_1\ldots\,k_n}[P]
=\dsize\msum{2}\Sb a_1,\,\ldots,\,a_\alpha\\ b_1,\,\ldots,\,b_\beta\endSb
\dsize\msum{2}\Sb \bar a_1,\,\ldots,\,\bar a_\nu\\ 
\bar b_1,\,\ldots,\,\bar b_\gamma\endSb
\dsize\msum{3}
\Sb c_1,\,\ldots,\,c_m\\ d_1,\,\ldots,\,d_n\endSb
\goth S^{\,i_1}_{a_1}\ldots\,\goth S^{\,i_\alpha}_{a_\alpha}\,\times\\
&\kern 40pt
\times\,\goth T^{b_1}_{j_1}\ldots\,\goth T^{b_\beta}_{j_\beta}
\ \overline{\goth S^{\,\bar i_1}_{\bar a_1}}
\ldots\,\overline{\goth S^{\,\bar i_\nu}_{\bar a_\nu}}\ 
\ \overline{\goth T^{\,\bar b_1}_{\bar j_1}}
\ldots\,\overline{\goth T^{\,\bar b_\gamma}_{\bar j_\gamma}}\ 
S^{h_1}_{c_1}\ldots\,S^{h_m}_{c_m}\,\times\\
\vspace{1.5ex}
&\kern 60pt
\times\,T^{\,d_1}_{k_1}\ldots\,T^{\,d_n}_{k_n}\ 
\tilde S^{\,a_1\ldots\,a_\alpha\,\bar a_1\ldots\,\bar a_\nu\,
c_1\ldots\,c_m}_{\,b_1\ldots\,b_\beta\,\bar b_1\ldots\,\bar b_\gamma
\,d_1\ldots\,d_n}[P],
\endaligned
\right.
\mytag{4.8}
\endalign
$$
where $\alpha=\alpha_P$, $\beta=\beta_P$, $\nu=\nu_P$, $\gamma
=\gamma_P$, $m=m_P$, $n=n_P$, and the integer number $P$ runs from 
$1$ to $J$. The components of transition matrices $\goth S$, $\goth T$,
$S=\varphi(\goth S)$, and $T=\varphi(\goth T)$ in \mythetag{4.7} and
\mythetag{4.8} are taken from the frame relationships \mythetag{1.3} 
and \mythetag{1.5}.
\head
5. Extended spinor connections.
\endhead
     Extended spinor connections were introduced in \mycite{1} in 
order to describe the structure of differentiations of extended 
spin-tensorial fields. In order to define them here we consider 
a pair of equipped local charts $U$ and $\tilde U$ with non-empty
intersection. Then we introduce the following $\theta$-parameters
defined within $U\cap\tilde U$:
$$
\align
&\hskip -6em
\tilde\theta^k_{ij}=\sum^3_{a=0}T^k_a\,L_{\tilde{\boldsymbol\Upsilon}_i}
\!(S^a_j)=\sum^3_{a=0}\sum^3_{v=0}T^k_a\,\tilde\Upsilon^v_i\,\frac{\partial
S^a_j}{\partial\tx^v}=-\sum^3_{a=0}L_{\tilde{\boldsymbol\Upsilon}_i}
\!(T^k_a)\,S^a_j,\hskip -2em
\mytag{5.1}\\
&\hskip -6em
\tilde\vartheta^k_{ij}=\sum^2_{a=1}\goth T^k_a\,
L_{\tilde{\boldsymbol\Upsilon}_i}(\goth S^a_j)
=\sum^2_{a=1}\sum^3_{v=0}\goth T^k_a\,\tilde\Upsilon^v_i\,
\frac{\partial\goth S^a_j}{\partial\tx^v}=
-\sum^2_{a=1}L_{\tilde{\boldsymbol\Upsilon}_i}
\!(\goth T^k_a)\,\goth S^a_j.\hskip -2em
\mytag{5.2}
\endalign
$$
The $\theta$-parameters without tilde are introduced in a similar way:
$$
\align
&\hskip -6em
\theta^k_{ij}=\sum^3_{a=0}S^k_a\,L_{\boldsymbol\Upsilon_i}
\!(T^a_j)=\sum^3_{a=0}\sum^3_{v=0}S^k_a\,\Upsilon^v_i\,\frac{\partial
T^a_j}{\partial x^v}=-\sum^3_{a=0}L_{\boldsymbol\Upsilon_i}
\!(S^k_a)\,T^a_j,\hskip -2em
\mytag{5.3}\\
&\hskip -6em
\vartheta^k_{ij}=\sum^2_{a=1}\goth S^k_a\,
L_{\boldsymbol\Upsilon_i}(\goth T^a_j)
=\sum^2_{a=1}\sum^3_{v=0}\goth S^k_a\,\Upsilon^v_i\,
\frac{\partial\goth T^a_j}{\partial\tx^v}=
-\sum^2_{a=1}L_{\boldsymbol\Upsilon_i}
\!(\goth S^k_a)\,\goth T^a_j.\hskip -2em
\mytag{5.4}
\endalign
$$
Note that $L_{\tilde{\boldsymbol\Upsilon}_i}$ and
$L_{\boldsymbol\Upsilon_i}$ in \mythetag{5.1}, \mythetag{5.2},
\mythetag{5.3}, and \mythetag{5.4} are Lie derivatives applied to
scalar functions $T^k_a$, $\goth T^k_a$, $S^k_a$, and $\goth S^k_a$
respectively. When applied to an arbitrary scalar function $f$
in $U\cap\tilde U$, the Lie derivative $L_{\boldsymbol\Upsilon_i}$ 
acts as follows:
$$
\hskip -2em
L_{\boldsymbol\Upsilon_i}(f)=\sum^3_{j=0}\Upsilon^j_i\,
\frac{\partial f}{\partial x^j}.
\mytag{5.5}
$$
The coefficients $\Upsilon^j_i$ in \mythetag{5.5} coincide with those
in the expansion \mythetag{1.13}.
\mydefinition{5.1} Let $N$ be the composite spin-tensorial bundle
\mythetag{4.2} over the space-time manifold $M$. An extended spinor
connection is a geometric object such that in each equipped local chart 
$U$ of $M$ it is represented by its components $\Alpha^k_{j\,i}$,
$\bar{\Alpha}\vphantom{\Alpha}^k_{j\,i}$, $\Gamma^k_{j\,i}$ and such 
that its components are smooth functions of the variables \mythetag{4.4} 
and \mythetag{4.5} transforming according to the formulas 
$$
\align
&\hskip -2em
\Alpha^k_{j\,i}=\dsize\sum^2_{b=1}\sum^2_{a=1}\sum^3_{c=0}
\goth S^k_a\,\goth T^b_i\,T^c_j\ \tilde{\Alpha}\vphantom{\Alpha}^a_{c\,b}
+\vartheta^k_{j\,i},
\mytag{5.6}\\
&\hskip -2em
\bar{\Alpha}\vphantom{\Alpha}^k_{j\,i}=\sum^2_{b=1}\sum^2_{a=1}
\sum^3_{c=0}\overline{\goth S^k_a}\ \overline{\goth T^b_i}\,T^c_j\  
\tilde{\bar{\Alpha}}\vphantom{\Alpha}^a_b
+\overline{\vartheta^k_{j\,i}},
\mytag{5.7}\\
\vspace{2ex}
&\hskip -2em
\Gamma^k_{j\,i}=\dsize\sum^3_{b=0}\sum^3_{a=0}\sum^3_{c=0}
S^k_a\,T^b_i\,T^c_j\ \tilde\Gamma^a_{c\,b}+\theta^k_{j\,i}
\mytag{5.8}
\endalign
$$
under a change of a local chart. The $\theta$-parameters in 
the transformation formulas \mythetag{5.6}, \mythetag{5.7}, 
and \mythetag{5.8} are taken from \mythetag{5.3} and
\mythetag{5.4}.
\enddefinition
    Extended spinor connections naturally arise when we describe
the set of differentiations of extended tensor fields. According
to the structural theorem proved in \mycite{1} each differentiation
$D$ is a sum of the three special types of differentiations:
\roster
\item a spatial covariant differentiation;
\item several native vertical multivariate differentiations;
\item a degenerate differentiation.
\endroster\par
    {\bf Native vertical multivariate differentiations} are most
simple in the above list in that sense that they are given by the
shortest formulas of all three. Let $\bold S[P]$ be a spin-tensor
from the list \mythetag{4.3} and let $(\alpha_P,\beta_P|\nu_P,
\gamma_P|m_P,n_P)$ be its type. Assume that $\bold Y$ is some
arbitrary extended spin-tensorial field of this type $(\alpha_P,
\beta_P|\nu_P,\gamma_P|m_P,n_P)$. Then for $\bold Y$ the native 
vertical multivariate differentiation $\vnabla_{\bold Y}[P]$ along 
this spin-tensorial field is defined (see \mycite{1}). In an equipped 
local chart $U$ it is represented by the native multivariate derivative
$$
\hskip -2em
\vnabla^{j_1\ldots\,j_\beta\bar j_1\ldots\,\bar j_\gamma
k_1\ldots\,k_n}_{i_1\ldots\,i_\alpha\bar i_1\ldots\,\bar i_\nu
h_1\ldots\,h_m}[P]=\frac{\partial}
{\partial S^{\,i_1\ldots\,i_\alpha\bar i_1\ldots\,\bar i_\nu h_1\ldots
\,h_m}_{j_1\ldots\,j_\beta\bar j_1\ldots\,\bar j_\gamma
k_1\ldots\,k_n}[P]
\vphantom{\vrule height 11pt depth 0pt}},
\mytag{5.9}
$$
where $\alpha=\alpha_P$, $\beta=\beta_P$, $\nu=\nu_P$, $\gamma=\gamma_P$,
$m=m_P$, $n=n_P$. Similarly, if $\bold Y$ is some extended spin-tensorial
field of the type $(\nu_P,\gamma_P|\alpha_P,\beta_P|m_P,n_P)$, then the
barred native vertical multivariate differentiation $\Bar\vnabla_{\bold
Y}[P]$ along $\bold Y$ is defined (see \mycite{1} again). In an equipped
local chart $U$ this differentiation is represented by the corresponding
barred native multivariate derivative
$$
\hskip -2em
\Bar\vnabla^{j_1\ldots\,j_\gamma\,\bar j_1\ldots\,\bar j_\beta\,
k_1\ldots\,k_n}_{i_1\ldots\,i_\nu\,\bar i_1\ldots\,\bar i_\alpha\,
h_1\ldots\,h_m}[P]=\frac{\partial}
{\partial\overline{S^{\,\bar i_1\ldots\,\bar i_\alpha\,i_1\ldots\,
i_\nu\,h_1\ldots\,h_m}_{\bar j_1\ldots\,\bar j_\beta\,j_1\ldots\,
j_\gamma\,k_1\ldots\,k_n}[P]
\vphantom{\vrule height 11pt depth 0pt}}}.
\mytag{5.10}
$$
Here again $\alpha=\alpha_P$, $\beta=\beta_P$, $\nu=\nu_P$, $\gamma
=\gamma_P$, $m=m_P$, $n=n_P$. The derivatives \mythetag{5.9} and
\mythetag{5.10} are called native because they are represented by 
partial derivatives with respect to the variables \mythetag{4.5} 
native for the composite bundle \mythetag{4.2}. Note that they
do not require and do not provide any geometric structures other
than those already exist due to the bundle \mythetag{4.2}.\par
     {\bf Degenerate differentiations} are given by 
a little bit more complicated formulas. Their structure is described
by the following theorem proved in \mycite{1}.
\mytheorem{5.1} Defining a degenerate differentiation $D$ of extended
spin-tensorial fields is equivalent to fixing three extended spin-tensorial
fields $\eufb S$, $\bar{\eufb S}$, and $\bold S$ of the types $(1,1|0,0|0,0)$,
$(0,0|1,1|0,0)$, and $(0,0|0,0|1,1)$ respectively.
\endproclaim
    Let $D$ be a degenerate differentiation and let $\bold X$ be an
arbitrary smooth spin-tensorial field of the type $(\varepsilon,\eta|
\sigma,\zeta|e,f)$. In an equipped local chart $U$ the components of
the field $D(\bold X)$ are given by the formula
$$
\gathered
DX^{a_1\ldots\,a_\varepsilon\bar a_1\ldots\,\bar a_\sigma c_1
\ldots\,c_e}_{b_1\ldots\,b_\eta\bar b_1\ldots\,\bar b_\zeta
d_1\ldots\,d_f}
=\sum^\varepsilon_{\mu=1}\sum^2_{v_\mu=1}\goth S^{a_\mu}_{v_\mu}\ 
X^{a_1\ldots\,v_\mu\,\ldots\,a_\varepsilon\bar a_1\ldots\,\bar a_\sigma
c_1\ldots\,c_e}_{b_1\ldots\,\ldots\,\ldots\,b_\eta\bar b_1\ldots\,
\bar b_\zeta d_1\ldots\,d_f}\,-\\
\kern 9em-\sum^\eta_{\mu=1}\sum^2_{w_\mu=1}\goth S^{w_\mu}_{b_\mu}\
X^{a_1\ldots\,\ldots\,\ldots\,a_\varepsilon\bar a_1\ldots\,\bar a_\sigma
c_1\ldots\,c_e}_{b_1\ldots\,w_\mu\,\ldots\,b_\eta\bar b_1\ldots\,
\bar b_\zeta d_1\ldots\,d_f}\,+\\
\kern -9em
+\sum^\sigma_{\mu=1}\sum^2_{v_\mu=1}
\bar{\goth S}\vphantom{\goth S}^{\bar a_\mu}_{v_\mu}\ 
X^{a_1\ldots\,a_\varepsilon\bar a_1\ldots\,v_\mu\,\ldots\,\bar a_\sigma
c_1\ldots\,c_e}_{b_1\ldots\,b_\eta\bar b_1\ldots\,\ldots\,\ldots\,
\bar b_\zeta d_1\ldots\,d_f}\,-\\
\kern 9em-\sum^\zeta_{\mu=1}\sum^2_{w_\mu=1}
\bar{\goth S}\vphantom{\goth S}^{w_\mu}_{\bar b_\mu}\
X^{a_1\ldots\,a_\varepsilon\bar a_1\ldots\,\ldots\,\ldots\,\bar a_\sigma
c_1\ldots\,c_e}_{b_1\ldots\,b_\eta\bar b_1\ldots\,w_\mu\,\ldots\,
\bar b_\zeta d_1\ldots\,d_f}\,+\\
\kern -9em+\sum^e_{\mu=1}\sum^3_{v_\mu=0}S^{c_\mu}_{v_\mu}\ 
X^{a_1\ldots\,a_\varepsilon\bar a_1\ldots\,\bar a_\sigma
c_1\ldots\,v_\mu\,\ldots\,c_e}_{b_1\ldots\,b_\eta\bar b_1\ldots\,
\bar b_\zeta d_1\ldots\,\ldots\,\ldots\,d_f}\,-\\
\kern 9em-\sum^f_{\mu=1}\sum^3_{w_\mu=0}S^{w_\mu}_{b_\mu}\
X^{a_1\ldots\,a_\varepsilon\bar a_1\ldots\,\bar a_\sigma
c_1\ldots\,\ldots\,\ldots\,c_e}_{b_1\ldots\,b_\eta\bar b_1\ldots\,
\bar b_\zeta d_1\ldots\,w_\mu\,\ldots\,d_f}.
\endgathered\qquad\quad
\mytag{5.11}
$$
Note that the formula \mythetag{5.11} has no derivatives at all. 
For this reason the differentiation $D$ given by this formula is
a degenerate differentiation. Note also that $\goth S^{a_\mu}_{v_\mu}$,
$\goth S^{w_\mu}_{b_\mu}$, $\bar{\goth S}\vphantom{\goth S}^{\bar
a_\mu}_{v_\mu}$, $\bar{\goth S}\vphantom{\goth S}^{w_\mu}_{\bar b_\mu}$,
$S^{c_\mu}_{v_\mu}$, and $S^{w_\mu}_{b_\mu}$ in \mythetag{5.11} are the
components of the extended spin-tensorial fields declared in the 
theorem~\mythetheorem{5.1}. Do not mix them with the components of
transition matrices $\goth S$ and $S$ taken from \mythetag{1.3} and 
\mythetag{1.5}.\par
    {\bf Spatial covariant differentiations} are most complicated of
the above three types of differentiations. In an equipped local chart
$U$ they are represented by the corresponding spatial covariant
derivatives:
$$
\allowdisplaybreaks
\gather
\gathered
\nabla_{\!\!j}X^{a_1\ldots\,a_\varepsilon\bar a_1\ldots\,\bar a_\sigma 
c_1\ldots\,c_e}_{b_1\ldots\,b_\eta\bar b_1\ldots\,\bar b_\zeta
d_1\ldots\,d_f}
=\sum^3_{k=0}\Upsilon^k_j\,
\frac{\partial X^{a_1\ldots\,a_\varepsilon\bar a_1\ldots\,\bar a_\sigma 
c_1\ldots\,c_e}_{b_1\ldots\,b_\eta\bar b_1\ldots\,\bar b_\zeta d_1\ldots\,
d_f}}{\partial x^k}\,-\\
\vspace{1.5ex}
-\sum^J_{P=1}\dsize\msum{2}\Sb i_1,\,\ldots,\,i_\alpha\\
j_1,\,\ldots,\,j_\beta\\ \bar i_1,\,\ldots,\,\bar i_\nu\\ 
\bar j_1,\,\ldots,\,\bar j_\gamma\endSb
\msum{3}\Sb h_1,\,\ldots,\,h_m\\ k_1,\,\ldots,\,k_n\endSb
\left(\,\shave{\sum^\alpha_{\mu=1}\sum^2_{v_\mu=1}}\Alpha^{i_\mu}_{j\,v_\mu}
\ S^{\,i_1\ldots\,v_\mu\,\ldots\,i_\alpha\,\bar i_1\ldots\,\bar i_\nu
\,h_1\ldots\,h_m}_{j_1\ldots\,\ldots\,\ldots\,j_\beta\,\bar j_1\ldots\,
\bar j_\gamma\,k_1\ldots\,k_n}[P]\,-\right.\\
\vspace{0.5ex plus 0.5ex minus 0.5ex}
-\sum^\beta_{\mu=1}\sum^2_{w_\mu=1}\Alpha^{w_\mu}_{j\,j_\mu}
\ S^{\,i_1\ldots\,\ldots\,\ldots\,i_\alpha\,\bar i_1\ldots\,\bar i_\nu
\,h_1\ldots\,h_m}_{j_1\ldots\,w_\mu\,\ldots\,j_\beta\,\bar j_1\ldots\,
\bar j_\gamma\,k_1\ldots\,k_n}[P]
+\sum^\nu_{\mu=1}\sum^2_{v_\mu=1}
\bar{\Alpha}\vphantom{\Alpha}^{\bar i_\mu}_{j\,v_\mu}\,\times\\
\vspace{0.5ex plus 0.5ex minus 0.5ex}
\times\,S^{\,i_1\ldots\,i_\alpha\,\bar i_1\ldots\,v_\mu\,\ldots\,
\bar i_\nu\,h_1\ldots\,h_m}_{j_1\ldots\,j_\beta\,\bar j_1\ldots\,\ldots
\,\ldots\,\bar j_\gamma\,k_1\ldots\,k_n}[P]-
\sum^\gamma_{\mu=1}\sum^2_{w_\mu=1}
\bar{\Alpha}\vphantom{\Alpha}^{w_\mu}_{j\,\bar j_\mu}\,\times\kern 4em\\
\vspace{0.5ex plus 0.5ex minus 0.5ex}
\times\,S^{i_1\ldots\,i_\alpha\,\bar i_1\ldots\,\ldots\,\ldots\,\bar i_\nu
\,h_1\ldots\,h_m}_{j_1\ldots\,j_\beta\,\bar j_1\ldots\,w_\mu\,\ldots\,
\bar j_\gamma\,k_1\ldots\,k_n}[P]
+\sum^m_{\mu=1}\sum^3_{v_\mu=0}\Gamma^{h_\mu}_{j\,v_\mu}\,\times\\
\vspace{0.5ex plus 0.5ex minus 0.5ex}
\times\,S^{i_1\ldots\,i_\alpha\,\bar i_1\ldots\,\bar i_\nu
\,h_1\ldots\,v_\mu\,\ldots\,h_m}_{j_1\ldots\,j_\beta\,\bar j_1\ldots\,
\bar j_\gamma\,k_1\ldots\,\ldots\,\ldots\,k_n}[P]-
\sum^n_{\mu=1}\sum^3_{w_\mu=0}\Gamma^{w_\mu}_{j\,k_\mu}\,\times\kern 4em\\
\vspace{0.5ex plus 0.5ex minus 0.5ex}
\left.\vphantom{\shave{\sum^\alpha_{\mu=1}\sum^2_{v_\mu=1}}}
\times\,S^{i_1\ldots\,i_\alpha\,\bar i_1\ldots\,\bar i_\nu\,
h_1\ldots\,\ldots\,\ldots\,h_m}_{j_1\ldots\,j_\beta\,\bar j_1\ldots\,
\bar j_\gamma\,k_1\ldots\,w_\mu\,\ldots\,k_n}[P]\right)
\frac{\partial X^{a_1\ldots\,a_\varepsilon\bar a_1\ldots\,\bar a_\sigma 
c_1\ldots\,c_e}_{b_1\ldots\,b_\eta\bar b_1\ldots\,\bar b_\zeta d_1\ldots\,
d_f}}{\partial S^{\,i_1\ldots\,i_\alpha\,\bar i_1\ldots\,\bar i_\nu\,
h_1\ldots\,h_m}_{j_1\ldots\,j_\beta\,\bar j_1\ldots\,\bar j_\gamma\,
k_1\ldots\,k_n}[P]}\,-\\
\vspace{0.5ex plus 0.5ex minus 0.5ex}
-\sum^J_{P=1}\dsize\msum{2}\Sb i_1,\,\ldots,\,i_\alpha\\
j_1,\,\ldots,\,j_\beta\\ \bar i_1,\,\ldots,\,\bar i_\nu\\ 
\bar j_1,\,\ldots,\,\bar j_\gamma\endSb
\msum{3}\Sb h_1,\,\ldots,\,h_m\\ k_1,\,\ldots,\,k_n\endSb
\left(\,\shave{\sum^\nu_{\mu=1}\sum^2_{v_\mu=1}}\Alpha^{i_\mu}_{j\,v_\mu}
\ \overline{S^{\,\bar i_1\ldots\,\bar i_\alpha\,i_1\ldots\,v_\mu\,\ldots\,
i_\nu h_1\ldots\,h_m}_{\bar j_1\ldots\,\bar j_\beta\,j_1\ldots\,\ldots\,
\ldots\,j_\gamma k_1\ldots\,k_n}[P]}\,-\right.\\
\vspace{0.5ex plus 0.5ex minus 0.5ex}
-\sum^\gamma_{\mu=1}\sum^2_{w_\mu=1}\Alpha^{w_\mu}_{j\,j_\mu}
\ \overline{S^{\,\bar i_1\ldots\,\bar i_\alpha\,i_1\ldots\,\ldots\,\ldots
\,i_\nu\,h_1\ldots\,h_m}_{\bar j_1\ldots\,\bar j_\beta\,j_1\ldots\,w_\mu
\,\ldots\,j_\gamma\,k_1\ldots\,k_n}[P]}+\sum^\alpha_{\mu=1}\sum^2_{v_\mu=1}
\bar{\Alpha}\vphantom{\Alpha}^{\bar i_\mu}_{j\,v_\mu}\,\times\\
\vspace{0.5ex plus 0.5ex minus 0.5ex}
\times\,\overline{S^{\,\bar i_1\ldots\,v_\mu\,\ldots\,\bar i_\alpha\,i_1
\ldots\,i_\nu\,h_1\ldots\,h_m}_{\bar j_1\ldots\,\ldots\,\ldots\,
\bar j_\beta\,j_1\ldots\,j_\gamma\,k_1\ldots\,k_n}[P]}
-\sum^\beta_{\mu=1}\sum^2_{w_\mu=1}
\bar{\Alpha}\vphantom{\Alpha}^{w_\mu}_{j\,\bar j_\mu}\,\times\\
\vspace{0.5ex plus 0.5ex minus 0.5ex}
\times\,\overline{S^{\,\bar i_1\ldots\,\ldots\,\ldots\,\bar i_\alpha\,
i_1\ldots\,i_\nu\,h_1\ldots\,h_m}_{\bar j_1\ldots\,w_\mu\,\ldots\,
\bar j_\beta\,j_1\ldots\,j_\gamma\,k_1\ldots\,k_n}[P]}+
\sum^m_{\mu=1}\sum^3_{v_\mu=0}\Gamma^{h_\mu}_{j\,v_\mu}\,\times
\endgathered\qquad
\mytag{5.12}\\
\gathered
\times\,\overline{S^{\bar i_1\ldots\,\bar i_\alpha\,i_1\ldots\,i_\nu\,
h_1\ldots\,v_\mu\,\ldots\,h_m}_{\bar j_1\ldots\,\bar j_\beta\,j_1
\ldots\,j_\gamma\,k_1\ldots\,\ldots\,\ldots\,k_n}[P]}-
\sum^n_{\mu=1}\sum^3_{w_\mu=0}\Gamma^{w_\mu}_{j\,k_\mu}\,\times\\
\left.\vphantom{\shave{\sum^\nu_{\mu=1}\sum^2_{v_\mu=1}}}
\times\,\overline{S^{\bar i_1\ldots\,\bar i_\alpha\,i_1\ldots\,i_\nu
\,h_1\ldots\,\ldots\,\ldots\,h_m}_{\bar j_1\ldots\,\bar j_\beta
j_1\ldots\,j_\gamma\,k_1\ldots\,w_\mu\,\ldots\,k_n}[P]}
\right)\frac{\partial X^{a_1\ldots\,a_\varepsilon\bar a_1\ldots\,
\bar a_\sigma c_1\ldots\,c_e}_{b_1\ldots\,b_\eta\bar b_1\ldots\,
\bar b_\zeta d_1\ldots\,d_f}}{\partial\overline{S^{\,i_1\ldots\,
i_\alpha\,\bar i_1\ldots\,\bar i_\nu\,h_1\ldots\,h_m}_{j_1\ldots\,
j_\beta\,\bar j_1\ldots\,\bar j_\gamma\,k_1\ldots\,k_n}[P]}\,}\,+\\
\kern -9em
+\sum^\varepsilon_{\mu=1}\sum^2_{v_\mu=1}\Alpha^{a_\mu}_{j\,v_\mu}\ 
X^{a_1\ldots\,v_\mu\,\ldots\,a_\varepsilon\bar a_1\ldots\,\bar a_\sigma
c_1\ldots\,c_e}_{b_1\ldots\,\ldots\,\ldots\,b_\eta\bar b_1\ldots\,
\bar b_\zeta d_1\ldots\,d_f}\,-\\
\kern 9em-\sum^\eta_{\mu=1}\sum^2_{w_\mu=1}\Alpha^{w_\mu}_{j\,b_\mu}\
X^{a_1\ldots\,\ldots\,\ldots\,a_\varepsilon\bar a_1\ldots\,\bar a_\sigma
c_1\ldots\,c_e}_{b_1\ldots\,w_\mu\,\ldots\,b_\eta\bar b_1\ldots\,
\bar b_\zeta d_1\ldots\,d_f}\,+\\
\kern -9em
+\sum^\sigma_{\mu=1}\sum^2_{v_\mu=1}
\bar{\Alpha}\vphantom{\Alpha}^{\bar a_\mu}_{j\,v_\mu}\ 
X^{a_1\ldots\,a_\varepsilon\bar a_1\ldots\,v_\mu\,\ldots\,\bar a_\sigma
c_1\ldots\,c_e}_{b_1\ldots\,b_\eta\bar b_1\ldots\,\ldots\,\ldots\,
\bar b_\zeta d_1\ldots\,d_f}\,-\\
\kern 9em-\sum^\zeta_{\mu=1}\sum^2_{w_\mu=1}
\bar{\Alpha}\vphantom{\Alpha}^{w_\mu}_{j\,\bar b_\mu}\
X^{a_1\ldots\,a_\varepsilon\bar a_1\ldots\,\ldots\,\ldots\,\bar a_\sigma
c_1\ldots\,c_e}_{b_1\ldots\,b_\eta\bar b_1\ldots\,w_\mu\,\ldots\,
\bar b_\zeta d_1\ldots\,d_f}\,+\\
\kern -9em+\sum^e_{\mu=1}\sum^3_{v_\mu=0}\Gamma^{c_\mu}_{j\,v_\mu}\ 
X^{a_1\ldots\,a_\varepsilon\bar a_1\ldots\,\bar a_\sigma
c_1\ldots\,v_\mu\,\ldots\,c_e}_{b_1\ldots\,b_\eta\bar b_1\ldots\,
\bar b_\zeta d_1\ldots\,\ldots\,\ldots\,d_f}\,-\\
\kern 9em-\sum^f_{\mu=1}\sum^3_{w_\mu=0}\Gamma^{w_\mu}_{j\,b_\mu}\
X^{a_1\ldots\,a_\varepsilon\bar a_1\ldots\,\bar a_\sigma
c_1\ldots\,\ldots\,\ldots\,c_e}_{b_1\ldots\,b_\eta\bar b_1\ldots\,
\bar b_\zeta d_1\ldots\,w_\mu\,\ldots\,d_f}.
\endgathered\kern 4em
\endgather
$$
The quantities $\Upsilon^k_j$ in the first term of \mythetag{5.12}
are taken from the expansion \mythetag{1.13}. The quantities
$\Alpha^k_{j\,i}$, $\bar{\Alpha}\vphantom{\Alpha}^k_{j\,i}$,
$\Gamma^k_{j\,i}$ in \mythetag{5.12} are the components of an
extended spinor connection introduced in the 
definition~\mythedefinition{5.1}.
\mytheorem{5.2} Defining a spatial covariant differentiation 
$\nabla$ of extended spin-tensorial fields is equivalent to 
defining some extended spinor connection.
\endproclaim
\noindent The theorem~\mythetheorem{5.2} is immediate from 
the formula \mythetag{5.12}. It was first proved in \mycite{1}.
\head
6. Commutation relationships and curvature spin-tensors.
\endhead
    According to the theorem~\mythetheorem{5.1} any degenerate
differentiation is given by three spin-tensorial fields $\eufb S$,
$\bar{\eufb S}$, and $\bold S$. Let's denote it as
$$
\hskip -2em
D=D(\eufb S,\bar{\eufb S},\bold S).
\mytag{6.1}
$$
Assume that we have two degenerate differentiations of the form
\mythetag{6.1}:
$$
\xalignat 2
&\hskip -2em
D_1=D(\eufb S_1,\bar{\eufb S}_1,\bold S_1),
&&D_2=D(\eufb S_2,\bar{\eufb S}_2,\bold S_2).
\mytag{6.2}
\endxalignat 
$$
The commutator of any two differentiations is a differentiation
(see section~12 in \mycite{1}). In the present case the commutator
of the degenerate differentiations $D_1$ and $D_2$ in \mythetag{6.2} 
is a degenerate differentiation:
$$
\hskip -2em
[D(\eufb S_1,\bar{\eufb S}_1,\bold S_1),\,D(\eufb S_2,\bar{\eufb
S}_2,\bold S_2)]=D(\eufb S_3,\bar{\eufb S}_3,\bold S_3).
\mytag{6.3}
$$
For the extended spin-tensorial fields $\eufb S_3$, $\bar{\eufb S}_3$, 
and $\bold S_3$ in \mythetag{6.3} one easily derives:
$$
\align
\eufb S_3&=C(\eufb S_1\otimes\eufb S_2-\eufb S_2\otimes\eufb S_1),\\
\vspace{1ex}
\bar{\eufb S}_3&=C(\bar{\eufb S}_1\otimes\bar{\eufb S}_2-\bar{\eufb S}_2
\otimes\bar{\eufb S}_1),\\
\vspace{1ex}
\bold S_3&=C(\bold S_1\otimes\bold S_2-\bold S_2\otimes\bold S_1).
\endalign
$$
These formulas mean that the operator-valued fields $\eufb S_3$, 
$\bar{\eufb S}_3$, and $\bold S_3$ are calculated as pointwise 
commutators of the corresponding fields $\eufb S_1$, $\bar{\eufb S}_1$,
$\bold S_1$ and $\eufb S_2$, $\bar{\eufb S}_2$, $\bold S_2$.\par
     Let $\vnabla_{\bold X}[P]$ be the $P$-th native multivariate
differentiation along an extended spin-tensorial field $\bold X$.
Then we have the equality
$$
\gather
\hskip -2em
[\vnabla_{\bold X}[P],\,D(\eufb S,\bar{\eufb S},\bold S)]=D(\eufb R,
\bar{\eufb R},\bold R),
\mytag{6.4}\\
\vspace{-1.5ex}
\intertext{where}
\vspace{-1.5ex}
\hskip -2em
\eufb R=\vnabla_{\bold X}[P]\eufb S,\qquad \bar{\eufb R}=
\vnabla_{\bold X}[P]\bar{\eufb S},\qquad \bold R=\vnabla_{\bold X}[P]
\bold S.
\mytag{6.5}
\endgather
$$
Similarly, for the $P$-th barred native multivariate
differentiation $\bar{\vnabla}_{\bold X}[P]$ along some extended
spin-tensorial field $\bold X$ we have the equality
$$
\gather
\hskip -2em
[\bar{\vnabla}_{\bold X}[P],\,D(\eufb S,\bar{\eufb S},\bold S)]
=D(\eufb R,\bar{\eufb R},\bold R),
\mytag{6.6}\\
\vspace{-1.5ex}
\intertext{where}
\vspace{-1.5ex}
\hskip -2em
\eufb R=\bar{\vnabla}_{\bold X}[P]\eufb S,\qquad \bar{\eufb R}=
\bar{\vnabla}_{\bold X}[P]\bar{\eufb S},\qquad \bold R
=\bar{\vnabla}_{\bold X}[P]\bold S.
\mytag{6.7}
\endgather
$$
The formulas \mythetag{6.4}, \mythetag{6.5}, \mythetag{6.6}, and
\mythetag{6.7} are proved by direct calculations using some equipped
local chart $U$.\par
    For mutual commutators of barred and non-barred native multivariate
differentiations one can easily derive the following formulas:
$$
\gather
\hskip -2em
\gathered
[\vnabla_{\bold X}[P],\,\vnabla_{\bold Y}[Q]]=\vnabla_{\bold V}[Q]-
\vnabla_{\bold U}[P],\\
\vspace{1.0ex}
\text{where}\qquad\bold V=\vnabla_{\bold X}[P]\bold Y\qquad\text{and}
\qquad\bold U=\vnabla_{\bold Y}[Q]\bold X;
\endgathered
\mytag{6.8}\\
\vspace{2.5ex}
\hskip -2em
\gathered
[\vnabla_{\bold X}[P],\,\bar{\vnabla}_{\bold Y}[Q]]
=\bar{\vnabla}_{\bold V}[Q]-\vnabla_{\bold U}[P],\\
\vspace{1.0ex}
\text{where}\qquad\bold V=\vnabla_{\bold X}[P]\bold Y\qquad\text{and}
\qquad\bold U=\bar{\vnabla}_{\bold Y}[Q]\bold X;
\endgathered
\mytag{6.9}\\
\vspace{2.5ex}
\hskip -2em
\gathered
[\bar{\vnabla}_{\bold X}[P],\,\bar{\vnabla}_{\bold Y}[Q]]
=\bar{\vnabla}_{\bold V}[Q]-\bar{\vnabla}_{\bold U}[P],\\
\vspace{1.0ex}
\text{where}\qquad\bold V=\bar{\vnabla}_{\bold X}[P]\bold Y\qquad
\text{and}\qquad\bold U=\bar{\vnabla}_{\bold Y}[Q]\bold X.
\endgathered
\mytag{6.10}
\endgather
$$
These formulas \mythetag{6.8}, \mythetag{6.9}, and \mythetag{6.10}
are also proved by direct calculations using some equipped local chart 
$U$ of $M$.\par 
     Now assume that we have some extended spinor connection associated
with the composite spin-tensorial bundle \mythetag{4.2}. \pagebreak
Then the spatial covariant differentiation $\nabla$ is defined and we can
write various commutation relationships with it:
$$
\gathered
[\nabla_{\bold X},\,D(\eufb S,\bar{\eufb S},\bold S)]
=D(\eufb R,\bar{\eufb R},\bold R),\\
\text{where}\qquad\eufb R=\nabla_{\bold X}\eufb S,
\qquad\bar{\eufb R}=\nabla_{\bold X}\bar{\eufb S},
\qquad \bold R=\nabla_{\bold X}\bold S.
\endgathered
\qquad
\mytag{6.11}
$$
The formulas for commutators of $\nabla_{\bold X}$ with barred and
non-barred multivariate differentiations are more complicated than
\mythetag{6.11}. In the case of $\vnabla_{\bold Y}[P]$ we have
$$
\gathered
[\nabla_{\!\bold X},\,\vnabla_{\bold Y}[P]]=
\vnabla_{\!\bold U}[P]+\sum^J_{Q=1}\vnabla_{\!\bold U[Q]}[Q]\,+\\
+\sum^J_{Q=1}\bar{\vnabla}_{\!\bar{\bold U}[Q]}[Q]
-\nabla_{\bold V}+D(\eufb N^{\sssize +},\bar{\eufb N}^{\sssize +},
\bold N^{\sssize +}),
\endgathered
\qquad
\mytag{6.12}
$$
where $\bold U=\nabla_{\bold X}\bold Y$ and $\bold V=\vnabla_{\bold Y}[P]
\bold X$. As for the fields $\bold U[Q]$ and $\bar{\bold U}[Q]$ in
\mythetag{6.12}, they are defined by the following formulas:
$$
\hskip -2em
\aligned
\bold U[Q]&=-D(\eufb N^{\sssize +},\bar{\eufb N}^{\sssize +},
\bold N^{\sssize +})\,\bold S[Q],\\
\vspace{1ex}
\bar{\bold U}[Q]&=-D(\eufb N^{\sssize +},\bar{\eufb N}^{\sssize +},
\bold N^{\sssize +})\,\tau(\bold S[Q]).
\endaligned
\mytag{6.13}
$$
Three spin-tensorial fields $\eufb N^{\sssize +}$, 
$\bar{\eufb N}^{\sssize +}$, $\bold N^{\sssize +}$ determining the
degenerate differentiation $D(\eufb N^{\sssize +},\bar{\eufb N}^{\sssize +},
\bold N^{\sssize +})$ in \mythetag{6.12} and \mythetag{6.13} are introduced
through three dynamic curvature spin-tensors $\eufb D^{\sssize +}[P]$, 
$\bar{\eufb D}^{\sssize +}[P]$, and $\bold D^{\sssize +}[P]$ respectively:
$$
\align
\hskip -2em
\eufb N^{\sssize +}&=\eufb D^{\sssize +}[P](\bold X,\bold Y)
=C(\eufb D^{\sssize +}[P]\otimes\bold X\otimes\bold Y),\\
\vspace{1ex}
\hskip -2em
\bar{\eufb N}^{\sssize +}&=\bar{\eufb D}^{\sssize +}[P](\bold X,
\bold Y)=C(\bar{\eufb D}^{\sssize +}[P]\otimes\bold X\otimes\bold Y),
\mytag{6.14}\\
\vspace{1ex}
\hskip -2em
\bold N^{\sssize +}&=\bold D^{\sssize +}[P](\bold X,\bold Y)
=C(\bold D^{\sssize +}[P]\otimes\bold X\otimes\bold Y).
\endalign
$$
The dynamic curvature spin-tensors $\eufb D^{\sssize +}[P]$, 
$\bar{\eufb D}^{\sssize +}[P]$, and $\bold D^{\sssize +}[P]$ in
\mythetag{6.14} are determined by the extended spinor connection 
through which the covariant differentiation $\nabla$ is defined. 
For the components of $\eufb D^{\sssize +}[P]$, 
$\bar{\eufb D}^{\sssize +}[P]$, and $\bold D^{\sssize +}[P]$ we 
get 
$$
\align
\hskip -2em
\goth D^{{\sssize +}\,k\,j_1\ldots\,j_\beta\bar j_1\ldots\,\bar j_\gamma
k_1\ldots\,k_n}_{i\,j\,i_1\ldots\,i_\alpha\bar i_1\ldots\,\bar i_\nu
h_1\ldots\,h_m}[P]&=-\frac{\partial\Alpha^k_{j\,i}}
{\partial S^{\,i_1\ldots\,i_\alpha\bar i_1\ldots\,\bar i_\nu h_1
\ldots\,h_m}_{j_1\ldots\,j_\beta\bar j_1\ldots\,\bar j_\gamma
k_1\ldots\,k_n}[P]\vphantom{\vrule height 11pt depth 0pt}},\\
\vspace{1ex}
\hskip -2em
\bar{\goth D}^{{\sssize +}\,k\,j_1\ldots\,j_\beta\bar j_1\ldots\,
\bar j_\gamma k_1\ldots\,k_n}_{i\,j\,i_1\ldots\,i_\alpha\bar i_1\ldots
\,\bar i_\nu h_1\ldots\,h_m}[P]&=-\frac{\partial\bar{\Alpha}
\vphantom{\Alpha}^k_{j\,i}}
{\partial S^{\,i_1\ldots\,i_\alpha\bar i_1\ldots\,\bar i_\nu h_1
\ldots\,h_m}_{j_1\ldots\,j_\beta\bar j_1\ldots\,\bar j_\gamma
k_1\ldots\,k_n}[P]\vphantom{\vrule height 11pt depth 0pt}},
\mytag{6.15}\\
\vspace{1ex}
\hskip -2em
D^{{\sssize +}\,k\,j_1\ldots\,j_\beta\bar j_1\ldots\,\bar j_\gamma
k_1\ldots\,k_n}_{i\,j\,i_1\ldots\,i_\alpha\bar i_1\ldots\,\bar i_\nu
h_1\ldots\,h_m}[P]&=-\frac{\partial\Gamma^k_{j\,i}}
{\partial S^{\,i_1\ldots\,i_\alpha\bar i_1\ldots\,\bar i_\nu h_1
\ldots\,h_m}_{j_1\ldots\,j_\beta\bar j_1\ldots\,\bar j_\gamma
k_1\ldots\,k_n}[P]\vphantom{\vrule height 11pt depth 0pt}}.
\endalign
$$\par
     In the case of the barred $P$-th native multivariate differentiation
along an extended spin-tensorial field $\bold Y$ the formula \mythetag{6.12}
is rewritten as follows:
$$
\gathered
[\nabla_{\!\bold X},\,\bar{\vnabla}_{\bold Y}[P]]=
\bar{\vnabla}_{\!\bold U}[P]+\sum^J_{Q=1}\vnabla_{\!\bold U[Q]}[Q]\,+\\
+\sum^J_{Q=1}\bar{\vnabla}_{\!\bar{\bold U}[Q]}[Q]
-\nabla_{\bold V}+D(\eufb N^{\sssize -},\bar{\eufb N}^{\sssize -},
\bold N^{\sssize -}).
\endgathered
\qquad
\mytag{6.16}
$$
Here $\bold U=\nabla_{\bold X}\bold Y$ and $\bold V
=\bar{\vnabla}_{\bold Y}[P]\bold X$. As for the fields $\bold U[Q]$ and
$\bar{\bold U}[Q]$ in \mythetag{6.16}, they are defined by the following
formulas:
$$
\hskip -2em
\aligned
\bold U[Q]&=-D(\eufb N^{\sssize -},\bar{\eufb N}^{\sssize -},
\bold N^{\sssize -})\,\bold S[Q],\\
\vspace{1ex}
\bar{\bold U}[Q]&=-D(\eufb N^{\sssize -},\bar{\eufb N}^{\sssize -},
\bold N^{\sssize -})\,\tau(\bold S[Q]).
\endaligned
\mytag{6.17}
$$
Three spin-tensorial fields $\eufb N^{\sssize -}$, 
$\bar{\eufb N}^{\sssize -}$, $\bold N^{\sssize -}$ determining the
degenerate differentiation $D(\eufb N^{\sssize -},\bar{\eufb N}^{\sssize -},
\bold N^{\sssize -})$ in \mythetag{6.16} and \mythetag{6.17} are introduced
through other three dynamic curvature spin-tensors 
$\eufb D^{\sssize -}[P]$, $\bar{\eufb D}^{\sssize -}[P]$, and 
$\bold D^{\sssize -}[P]$:
$$
\align
\hskip -2em
\eufb N^{\sssize -}&=\eufb D^{\sssize -}[P](\bold X,\bold Y)
=C(\eufb D^{\sssize -}[P]\otimes\bold X\otimes\bold Y),\\
\vspace{1ex}
\hskip -2em
\bar{\eufb N}^{\sssize -}&=\bar{\eufb D}^{\sssize -}[P](\bold X,
\bold Y)=C(\bar{\eufb D}^{\sssize -}[P]\otimes\bold X\otimes\bold Y),
\mytag{6.18}\\
\vspace{1ex}
\hskip -2em
\bold N^{\sssize -}&=\bold D^{\sssize -}[P](\bold X,\bold Y)
=C(\bold D^{\sssize -}[P]\otimes\bold X\otimes\bold Y).
\endalign
$$
Like $\eufb D^{\sssize +}[P]$, $\bar{\eufb D}^{\sssize +}[P]$, and 
$\bold D^{\sssize +}[P]$ in \mythetag{6.14}, the dynamic curvature
spin-tensors $\eufb D^{\sssize -}[P]$, $\bar{\eufb D}^{\sssize -}[P]$, 
and $\bold D^{\sssize -}[P]$ in \mythetag{6.18} are determined by the
extended spinor connection through which the covariant differentiation
$\nabla$ is defined. For the components of these spin-tensors we get
the following expressions:
$$
\align
\hskip -2em
\goth D^{{\sssize -}\,k\,j_1\ldots\,j_\gamma\,\bar j_1\ldots\,
\bar j_\beta\,k_1\ldots\,k_n}_{i\,j\,i_1\ldots\,i_\nu\,\bar i_1
\ldots\,\bar i_\alpha\,h_1\ldots\,h_m}[P]&=
-\frac{\partial\Alpha^k_{j\,i}}{\partial\overline{S^{\,\bar i_1
\ldots\,\bar i_\alpha\,i_1\ldots\,i_\nu\,h_1\ldots\,h_m}_{\bar j_1
\ldots\,\bar j_\beta\,j_1\ldots\,j_\gamma\,k_1\ldots
\,k_n}[P]\vphantom{\vrule height 11pt depth 0pt}}},\\
\vspace{1ex}
\hskip -2em
\bar{\goth D}^{{\sssize -}\,k\,j_1\ldots\,j_\gamma\,\bar j_1\ldots\,
\bar j_\beta\,k_1\ldots\,k_n}_{i\,j\,i_1\ldots\,i_\nu\,\bar i_1
\ldots\,\bar i_\alpha\,h_1\ldots\,h_m}[P]&=
-\frac{\partial\bar{\Alpha}\vphantom{\Alpha}^k_{j\,i}}
{\partial\overline{S^{\,\bar i_1\ldots\,\bar i_\alpha\,i_1\ldots
\,i_\nu\,h_1\ldots\,h_m}_{\bar j_1\ldots\,\bar j_\beta\,j_1\ldots
\,j_\gamma\,k_1\ldots\,k_n}[P]\vphantom{\vrule height 11pt depth 0pt}}},
\mytag{6.19}\\
\vspace{1ex}
\hskip -2em
D^{{\sssize -}\,k\,j_1\ldots\,j_\gamma\,\bar j_1\ldots\,
\bar j_\beta\,k_1\ldots\,k_n}_{i\,j\,i_1\ldots\,i_\nu\,\bar i_1
\ldots\,\bar i_\alpha\,h_1\ldots\,h_m}[P]&=
-\frac{\partial\Gamma^k_{j\,i}}{\partial\overline{S^{\,\bar i_1
\ldots\,\bar i_\alpha\,i_1\ldots\,i_\nu\,h_1\ldots\,h_m}_{\bar j_1
\ldots\,\bar j_\beta\,j_1\ldots\,j_\gamma\,k_1\ldots
\,k_n}[P]\vphantom{\vrule height 11pt depth 0pt}}}.
\endalign
$$
Using \mythetag{6.15} and \mythetag{6.19}, we can write explicit 
expressions for the components of the spin-tensorial fields
$\eufb N^{\sssize +}$, $\bar{\eufb N}^{\sssize +}$, $\bold N^{\sssize +}$,
$\eufb N^{\sssize -}$, $\bar{\eufb N}^{\sssize -}$, and
$\bold N^{\sssize -}$:
$$
\gather
\goth N^{{\sssize +}\,k}_{\ i}=\sum^3_{j=1}
\dsize\msums{2}{3}\Sb i_1,\,\ldots,\,i_\alpha\\
j_1,\,\ldots,\,j_\beta\\ \bar i_1,\,\ldots,\,\bar i_\nu\\ 
\bar j_1,\,\ldots,\,\bar j_\gamma\\
h_1,\,\ldots,\,h_m\\ k_1,\,\ldots,\,k_n\endSb
\goth D^{{\sssize +}\,k\,j_1\ldots\,j_\beta\bar j_1\ldots\,
\bar j_\gamma k_1\ldots\,k_n}_{i\,j\,i_1\ldots\,i_\alpha
\bar i_1\ldots\,\bar i_\nu h_1\ldots\,h_m}[P]\ X^j\ 
Y^{\,i_1\ldots\,i_\alpha\bar i_1\ldots\,\bar i_\nu h_1\ldots
\,h_m}_{j_1\ldots\,j_\beta\bar j_1\ldots\,\bar j_\gamma
k_1\ldots\,k_n},\\
\goth N^{{\sssize -}\,k}_{\ i}=\sum^3_{j=1}
\dsize\msums{2}{3}\Sb i_1,\,\ldots,\,i_\alpha\\
j_1,\,\ldots,\,j_\beta\\ \bar i_1,\,\ldots,\,\bar i_\nu\\ 
\bar j_1,\,\ldots,\,\bar j_\gamma\\
h_1,\,\ldots,\,h_m\\ k_1,\,\ldots,\,k_n\endSb
\goth D^{{\sssize -}\,k\,j_1\ldots\,j_\gamma\,\bar j_1\ldots\,
\bar j_\beta\,k_1\ldots\,k_n}_{i\,j\,i_1\ldots\,i_\nu\,\bar i_1
\ldots\,\bar i_\alpha\,h_1\ldots\,h_m}[P]\ X^j\ 
Y^{i_1\ldots\,i_\nu\,\bar i_1\ldots\,\bar i_\alpha\,
h_1\ldots\,h_m}_{j_1\ldots\,j_\gamma\,\bar j_1\ldots\,
\bar j_\beta\,k_1\ldots\,k_n}.
\endgather
$$
These two formulas are coordinate representations for the first
two formulas in \mythetag{6.14} and \mythetag{6.18}. Coordinate 
representations for other formulas in \mythetag{6.14} and 
\mythetag{6.18} are easily written by analogy.\par
     Again assume that we have some extended spinor connection 
associated with the composite spin-tensorial bundle \mythetag{4.2}
and suppose that $\nabla$ is the spatial covariant differentiation 
defined with the use this extended connection. Then
$$
[\nabla_{\bold X},\,\nabla_{\bold Y}]=
\nabla_{\!\bold U}+\sum^J_{Q=1}\vnabla_{\!\bold U[Q]}[Q]
+\sum^J_{Q=1}\bar{\vnabla}_{\!\bar{\bold U}[Q]}[Q]
+D(\eufb N,\bar{\eufb N},\bold N),
\quad
\mytag{6.20}
$$
where $\bold U=\nabla_{\bold X}\bold Y-\nabla_{\bold Y}\bold X
-\bold T(\bold X,\bold Y)$ and 
$$
\hskip -2em
\bold T(\bold X,\bold Y)=C(\bold T\otimes\bold X\otimes\bold Y).
\mytag{6.21}
$$
The extended spin-tensorial field $\bold T$ of the type $(0,0|0,0|1,2)$
in \mythetag{6.21} is known as the {\it torsion field}. The components
of the torsion are given by the formula
$$
\hskip -2em
T^k_{ij}=\Gamma^k_{ij}-\Gamma^k_{j\,i}-c^{\,k}_{ij},
\mytag{6.22}
$$
where $c^{\,k}_{ij}$ are taken from \mythetag{1.14} or from \mythetag{1.15}.
The extended spin tensorial fields $\bold U[Q]$ and $\bar{\bold U}[Q]$
in \mythetag{6.20} are obtained by applying $D(\eufb N,\bar{\eufb N},
\bold N)$ to the native fields:
$$
\hskip -2em
\aligned
\bold U[Q]&=-D(\eufb N,\bar{\eufb N},\bold N)\,\bold S[Q],\\
\vspace{1ex}
\bar{\bold U}[Q]&=-D(\eufb N,\bar{\eufb N},\bold N)\,\tau(\bold S[Q]).
\endaligned
\mytag{6.23}
$$
The degenerate differentiation $D(\eufb N,\bar{\eufb N},\bold N)$ in
\mythetag{6.20} and \mythetag{6.23} is determined by three 
spin-tensorial fields $\eufb N$, $\bar{\eufb N}$, $\bold N$ which are
expressed through three non-dynamic curvature spin-tensors
$\eufb R$, $\bar{\eufb R}$, and $\bold R$ respectively:
$$
\align
\hskip -2em
\eufb N&=\eufb R(\bold X,\bold Y)
=C(\eufb R\otimes\bold X\otimes\bold Y),\\
\vspace{1ex}
\hskip -2em
\bar{\eufb N}&=\bar{\eufb R}(\bold X,
\bold Y)=C(\bar{\eufb R}\otimes\bold X\otimes\bold Y),
\mytag{6.24}\\
\vspace{1ex}
\hskip -2em
\bold N&=\bold R(\bold X,\bold Y)
=C(\bold R\otimes\bold X\otimes\bold Y).
\endalign
$$
The components of the spin-curvature tensor $\eufb R$ in \mythetag{6.24} 
are given by the formula
$$
\allowdisplaybreaks
\gather
\gathered
\goth R^p_{qij}=\sum^3_{k=0}\Upsilon^k_i\,\frac{\partial\Alpha^p_{j\,q}}
{\partial x^k}-\sum^3_{k=0}\Upsilon^k_j\,\frac{\partial\Alpha^p_{i\,q}}
{\partial x^k}+\sum^2_{h=1}\left(\Alpha^p_{i\,h}\,\Alpha^{\!h}_{j\,q}
-\Alpha^p_{j\,h}\,\Alpha^{\!h}_{i\,q}\right)\,-\\
-\sum^J_{P=1}\dsize\msum{2}\Sb i_1,\,\ldots,\,i_\alpha\\
j_1,\,\ldots,\,j_\beta\\ \bar i_1,\,\ldots,\,\bar i_\nu\\ 
\bar j_1,\,\ldots,\,\bar j_\gamma\endSb
\msum{3}\Sb h_1,\,\ldots,\,h_m\\ k_1,\,\ldots,\,k_n\endSb
\left(\,\shave{\sum^\alpha_{\mu=1}\sum^2_{v_\mu=1}}\Alpha^{i_\mu}_{i\,v_\mu}
\ S^{\,i_1\ldots\,v_\mu\,\ldots\,i_\alpha\,\bar i_1\ldots\,\bar i_\nu
\,h_1\ldots\,h_m}_{j_1\ldots\,\ldots\,\ldots\,j_\beta\,\bar j_1\ldots\,
\bar j_\gamma\,k_1\ldots\,k_n}[P]\,-\right.\\
\vspace{0.5ex plus 0.5ex minus 0.5ex}
-\sum^\beta_{\mu=1}\sum^2_{w_\mu=1}\Alpha^{w_\mu}_{i\,j_\mu}
\ S^{\,i_1\ldots\,\ldots\,\ldots\,i_\alpha\,\bar i_1\ldots\,\bar i_\nu
\,h_1\ldots\,h_m}_{j_1\ldots\,w_\mu\,\ldots\,j_\beta\,\bar j_1\ldots\,
\bar j_\gamma\,k_1\ldots\,k_n}[P]
+\sum^\nu_{\mu=1}\sum^2_{v_\mu=1}
\bar{\Alpha}\vphantom{\Alpha}^{\bar i_\mu}_{i\,v_\mu}\,\times\\
\vspace{0.5ex plus 0.5ex minus 0.5ex}
\times\,S^{\,i_1\ldots\,i_\alpha\,\bar i_1\ldots\,v_\mu\,\ldots\,
\bar i_\nu\,h_1\ldots\,h_m}_{j_1\ldots\,j_\beta\,\bar j_1\ldots\,\ldots
\,\ldots\,\bar j_\gamma\,k_1\ldots\,k_n}[P]-
\sum^\gamma_{\mu=1}\sum^2_{w_\mu=1}
\bar{\Alpha}\vphantom{\Alpha}^{w_\mu}_{i\,\bar j_\mu}\,\times\\
\vspace{0.5ex plus 0.5ex minus 0.5ex}
\endgathered\qquad
\mytag{6.25}\\
\gathered
\times\,S^{i_1\ldots\,i_\alpha\,\bar i_1\ldots\,\ldots\,\ldots\,\bar i_\nu
\,h_1\ldots\,h_m}_{j_1\ldots\,j_\beta\,\bar j_1\ldots\,w_\mu\,\ldots\,
\bar j_\gamma\,k_1\ldots\,k_n}[P]
+\sum^m_{\mu=1}\sum^3_{v_\mu=0}\Gamma^{h_\mu}_{i\,v_\mu}\,\times\kern 4em\\
\vspace{0.5ex plus 0.5ex minus 0.5ex}
\times\,S^{i_1\ldots\,i_\alpha\,\bar i_1\ldots\,\bar i_\nu
\,h_1\ldots\,v_\mu\,\ldots\,h_m}_{j_1\ldots\,j_\beta\,\bar j_1\ldots\,
\bar j_\gamma\,k_1\ldots\,\ldots\,\ldots\,k_n}[P]-
\sum^n_{\mu=1}\sum^3_{w_\mu=0}\Gamma^{w_\mu}_{i\,k_\mu}\,\times\\
\vspace{0.5ex plus 0.5ex minus 0.5ex}
\left.\vphantom{\shave{\sum^\alpha_{\mu=1}\sum^2_{v_\mu=1}}}
\times\,S^{i_1\ldots\,i_\alpha\,\bar i_1\ldots\,\bar i_\nu\,
h_1\ldots\,\ldots\,\ldots\,h_m}_{j_1\ldots\,j_\beta\,\bar j_1\ldots\,
\bar j_\gamma\,k_1\ldots\,w_\mu\,\ldots\,k_n}[P]\right)
\frac{\partial\Alpha^p_{j\,q}}{\partial S^{\,i_1\ldots\,i_\alpha\,
\bar i_1\ldots\,\bar i_\nu\,
h_1\ldots\,h_m}_{j_1\ldots\,j_\beta\,\bar j_1\ldots\,\bar j_\gamma\,
k_1\ldots\,k_n}[P]}\,-\\
\vspace{0.5ex plus 0.5ex minus 0.5ex}
-\sum^J_{P=1}\dsize\msum{2}\Sb i_1,\,\ldots,\,i_\alpha\\
j_1,\,\ldots,\,j_\beta\\ \bar i_1,\,\ldots,\,\bar i_\nu\\ 
\bar j_1,\,\ldots,\,\bar j_\gamma\endSb
\msum{3}\Sb h_1,\,\ldots,\,h_m\\ k_1,\,\ldots,\,k_n\endSb
\left(\,\shave{\sum^\nu_{\mu=1}\sum^2_{v_\mu=1}}\Alpha^{i_\mu}_{i\,v_\mu}
\ \overline{S^{\,\bar i_1\ldots\,\bar i_\alpha\,i_1\ldots\,v_\mu\,\ldots\,
i_\nu h_1\ldots\,h_m}_{\bar j_1\ldots\,\bar j_\beta\,j_1\ldots\,\ldots\,
\ldots\,j_\gamma k_1\ldots\,k_n}[P]}\,-\right.\\
\vspace{0.5ex plus 0.5ex minus 0.5ex}
-\sum^\gamma_{\mu=1}\sum^2_{w_\mu=1}\Alpha^{w_\mu}_{i\,j_\mu}
\ \overline{S^{\,\bar i_1\ldots\,\bar i_\alpha\,i_1\ldots\,\ldots\,\ldots
\,i_\nu\,h_1\ldots\,h_m}_{\bar j_1\ldots\,\bar j_\beta\,j_1\ldots\,w_\mu
\,\ldots\,j_\gamma\,k_1\ldots\,k_n}[P]}+\sum^\alpha_{\mu=1}\sum^2_{v_\mu=1}
\bar{\Alpha}\vphantom{\Alpha}^{\bar i_\mu}_{i\,v_\mu}\,\times\\
\vspace{0.5ex plus 0.5ex minus 0.5ex}
\times\,\overline{S^{\,\bar i_1\ldots\,v_\mu\,\ldots\,\bar i_\alpha\,i_1
\ldots\,i_\nu\,h_1\ldots\,h_m}_{\bar j_1\ldots\,\ldots\,\ldots\,
\bar j_\beta\,j_1\ldots\,j_\gamma\,k_1\ldots\,k_n}[P]}
-\sum^\beta_{\mu=1}\sum^2_{w_\mu=1}
\bar{\Alpha}\vphantom{\Alpha}^{w_\mu}_{i\,\bar j_\mu}\,\times\\
\vspace{0.5ex plus 0.5ex minus 0.5ex}
\times\,\overline{S^{\,\bar i_1\ldots\,\ldots\,\ldots\,\bar i_\alpha\,
i_1\ldots\,i_\nu\,h_1\ldots\,h_m}_{\bar j_1\ldots\,w_\mu\,\ldots\,
\bar j_\beta\,j_1\ldots\,j_\gamma\,k_1\ldots\,k_n}[P]}+
\sum^m_{\mu=1}\sum^3_{v_\mu=0}\Gamma^{h_\mu}_{i\,v_\mu}\,\times\kern 4em\\
\times\,\overline{S^{\bar i_1\ldots\,\bar i_\alpha\,i_1\ldots\,i_\nu\,
h_1\ldots\,v_\mu\,\ldots\,h_m}_{\bar j_1\ldots\,\bar j_\beta\,j_1
\ldots\,j_\gamma\,k_1\ldots\,\ldots\,\ldots\,k_n}[P]}-
\sum^n_{\mu=1}\sum^3_{w_\mu=0}\Gamma^{w_\mu}_{i\,k_\mu}\,\times\\
\vspace{0.5ex plus 0.5ex minus 0.5ex}
\left.\vphantom{\shave{\sum^\nu_{\mu=1}\sum^2_{v_\mu=1}}}
\times\,\overline{S^{\bar i_1\ldots\,\bar i_\alpha\,i_1\ldots\,i_\nu
\,h_1\ldots\,\ldots\,\ldots\,h_m}_{\bar j_1\ldots\,\bar j_\beta
j_1\ldots\,j_\gamma\,k_1\ldots\,w_\mu\,\ldots\,k_n}[P]}
\right)\frac{\partial\Alpha^p_{j\,q}}{\partial\overline{S^{\,i_1\ldots
\,i_\alpha\,\bar i_1\ldots\,\bar i_\nu\,h_1\ldots\,h_m}_{j_1\ldots
\,j_\beta\,\bar j_1\ldots\,\bar j_\gamma\,k_1\ldots\,k_n}[P]}\,}\,+\\
\endgathered\kern 4em\\
\vspace{1ex plus 0.5ex minus 0.5ex}
\gathered
+\sum^J_{P=1}\dsize\msum{2}\Sb i_1,\,\ldots,\,i_\alpha\\
j_1,\,\ldots,\,j_\beta\\ \bar i_1,\,\ldots,\,\bar i_\nu\\ 
\bar j_1,\,\ldots,\,\bar j_\gamma\endSb
\msum{3}\Sb h_1,\,\ldots,\,h_m\\ k_1,\,\ldots,\,k_n\endSb
\left(\,\shave{\sum^\alpha_{\mu=1}\sum^2_{v_\mu=1}}\Alpha^{i_\mu}_{j\,v_\mu}
\ S^{\,i_1\ldots\,v_\mu\,\ldots\,i_\alpha\,\bar i_1\ldots\,\bar i_\nu
\,h_1\ldots\,h_m}_{j_1\ldots\,\ldots\,\ldots\,j_\beta\,\bar j_1\ldots\,
\bar j_\gamma\,k_1\ldots\,k_n}[P]\,-\right.\\
\vspace{0.5ex plus 0.5ex minus 0.5ex}
-\sum^\beta_{\mu=1}\sum^2_{w_\mu=1}\Alpha^{w_\mu}_{j\,j_\mu}
\ S^{\,i_1\ldots\,\ldots\,\ldots\,i_\alpha\,\bar i_1\ldots\,\bar i_\nu
\,h_1\ldots\,h_m}_{j_1\ldots\,w_\mu\,\ldots\,j_\beta\,\bar j_1\ldots\,
\bar j_\gamma\,k_1\ldots\,k_n}[P]
+\sum^\nu_{\mu=1}\sum^2_{v_\mu=1}
\bar{\Alpha}\vphantom{\Alpha}^{\bar i_\mu}_{j\,v_\mu}\,\times\\
\vspace{0.5ex plus 0.5ex minus 0.5ex}
\times\,S^{\,i_1\ldots\,i_\alpha\,\bar i_1\ldots\,v_\mu\,\ldots\,
\bar i_\nu\,h_1\ldots\,h_m}_{j_1\ldots\,j_\beta\,\bar j_1\ldots\,\ldots
\,\ldots\,\bar j_\gamma\,k_1\ldots\,k_n}[P]-
\sum^\gamma_{\mu=1}\sum^2_{w_\mu=1}
\bar{\Alpha}\vphantom{\Alpha}^{w_\mu}_{j\,\bar j_\mu}\,\times\\
\vspace{0.5ex plus 0.5ex minus 0.5ex}
\times\,S^{i_1\ldots\,i_\alpha\,\bar i_1\ldots\,\ldots\,\ldots\,\bar i_\nu
\,h_1\ldots\,h_m}_{j_1\ldots\,j_\beta\,\bar j_1\ldots\,w_\mu\,\ldots\,
\bar j_\gamma\,k_1\ldots\,k_n}[P]
+\sum^m_{\mu=1}\sum^3_{v_\mu=0}\Gamma^{h_\mu}_{j\,v_\mu}\,\times\kern 4em\\
\vspace{0.5ex plus 0.5ex minus 0.5ex}
\times\,S^{i_1\ldots\,i_\alpha\,\bar i_1\ldots\,\bar i_\nu
\,h_1\ldots\,v_\mu\,\ldots\,h_m}_{j_1\ldots\,j_\beta\,\bar j_1\ldots\,
\bar j_\gamma\,k_1\ldots\,\ldots\,\ldots\,k_n}[P]-
\sum^n_{\mu=1}\sum^3_{w_\mu=0}\Gamma^{w_\mu}_{j\,k_\mu}\,\times
\endgathered\kern 4em\\
\gathered
%\vspace{0.5ex plus 0.5ex minus 0.5ex}
\left.\vphantom{\shave{\sum^\alpha_{\mu=1}\sum^2_{v_\mu=1}}}
\times\,S^{i_1\ldots\,i_\alpha\,\bar i_1\ldots\,\bar i_\nu\,
h_1\ldots\,\ldots\,\ldots\,h_m}_{j_1\ldots\,j_\beta\,\bar j_1\ldots\,
\bar j_\gamma\,k_1\ldots\,w_\mu\,\ldots\,k_n}[P]\right)
\frac{\partial\Alpha^p_{i\,q}}{\partial S^{\,i_1\ldots\,i_\alpha\,
\bar i_1\ldots\,\bar i_\nu\,
h_1\ldots\,h_m}_{j_1\ldots\,j_\beta\,\bar j_1\ldots\,\bar j_\gamma\,
k_1\ldots\,k_n}[P]}\,+\\
+\sum^J_{P=1}\dsize\msum{2}\Sb i_1,\,\ldots,\,i_\alpha\\
j_1,\,\ldots,\,j_\beta\\ \bar i_1,\,\ldots,\,\bar i_\nu\\ 
\bar j_1,\,\ldots,\,\bar j_\gamma\endSb
\msum{3}\Sb h_1,\,\ldots,\,h_m\\ k_1,\,\ldots,\,k_n\endSb
\left(\,\shave{\sum^\nu_{\mu=1}\sum^2_{v_\mu=1}}\Alpha^{i_\mu}_{j\,v_\mu}
\ \overline{S^{\,\bar i_1\ldots\,\bar i_\alpha\,i_1\ldots\,v_\mu\,\ldots\,
i_\nu h_1\ldots\,h_m}_{\bar j_1\ldots\,\bar j_\beta\,j_1\ldots\,\ldots\,
\ldots\,j_\gamma k_1\ldots\,k_n}[P]}\,-\right.\\
\vspace{0.5ex plus 0.5ex minus 0.5ex}
-\sum^\gamma_{\mu=1}\sum^2_{w_\mu=1}\Alpha^{w_\mu}_{j\,j_\mu}
\ \overline{S^{\,\bar i_1\ldots\,\bar i_\alpha\,i_1\ldots\,\ldots\,\ldots
\,i_\nu\,h_1\ldots\,h_m}_{\bar j_1\ldots\,\bar j_\beta\,j_1\ldots\,w_\mu
\,\ldots\,j_\gamma\,k_1\ldots\,k_n}[P]}+\sum^\alpha_{\mu=1}\sum^2_{v_\mu=1}
\bar{\Alpha}\vphantom{\Alpha}^{\bar i_\mu}_{j\,v_\mu}\,\times\\
\vspace{0.5ex plus 0.5ex minus 0.5ex}
\times\,\overline{S^{\,\bar i_1\ldots\,v_\mu\,\ldots\,\bar i_\alpha\,i_1
\ldots\,i_\nu\,h_1\ldots\,h_m}_{\bar j_1\ldots\,\ldots\,\ldots\,
\bar j_\beta\,j_1\ldots\,j_\gamma\,k_1\ldots\,k_n}[P]}
-\sum^\beta_{\mu=1}\sum^2_{w_\mu=1}
\bar{\Alpha}\vphantom{\Alpha}^{w_\mu}_{j\,\bar j_\mu}\,\times\\
\vspace{0.5ex plus 0.5ex minus 0.5ex}
\times\,\overline{S^{\,\bar i_1\ldots\,\ldots\,\ldots\,\bar i_\alpha\,
i_1\ldots\,i_\nu\,h_1\ldots\,h_m}_{\bar j_1\ldots\,w_\mu\,\ldots\,
\bar j_\beta\,j_1\ldots\,j_\gamma\,k_1\ldots\,k_n}[P]}+
\sum^m_{\mu=1}\sum^3_{v_\mu=0}\Gamma^{h_\mu}_{j\,v_\mu}\,\times\kern 4em\\
\times\,\overline{S^{\bar i_1\ldots\,\bar i_\alpha\,i_1\ldots\,i_\nu\,
h_1\ldots\,v_\mu\,\ldots\,h_m}_{\bar j_1\ldots\,\bar j_\beta\,j_1
\ldots\,j_\gamma\,k_1\ldots\,\ldots\,\ldots\,k_n}[P]}-
\sum^n_{\mu=1}\sum^3_{w_\mu=0}\Gamma^{w_\mu}_{j\,k_\mu}\,\times\\
\left.\vphantom{\shave{\sum^\nu_{\mu=1}\sum^2_{v_\mu=1}}}
\times\,\overline{S^{\bar i_1\ldots\,\bar i_\alpha\,i_1\ldots\,i_\nu
\,h_1\ldots\,\ldots\,\ldots\,h_m}_{\bar j_1\ldots\,\bar j_\beta
j_1\ldots\,j_\gamma\,k_1\ldots\,w_\mu\,\ldots\,k_n}[P]}
\right)\frac{\partial\Alpha^p_{i\,q}}{\partial\overline{S^{\,i_1\ldots
\,i_\alpha\,\bar i_1\ldots\,\bar i_\nu\,h_1\ldots\,h_m}_{j_1\ldots
\,j_\beta\,\bar j_1\ldots\,\bar j_\gamma\,k_1\ldots\,k_n}[P]}\,}
-\sum^3_{k=0}c^{\,k}_{ij}\,\Alpha^p_{kq}.
\endgathered\kern 4em
\endgather
$$
For each $P$ in the above formula \mythetag{6.25} we implicitly assume
that $\alpha=\alpha_P$, $\beta=\beta_P$, $\nu=\nu_P$, $\gamma=\gamma_P$,
$m=m_P$, $n=n_P$. The components of $\bar{\eufb R}$ are calculated
similarly:
%\vskip -1ex
$$
\allowdisplaybreaks
\gather
\gathered
\bar{\goth R}^p_{qij}=\sum^3_{k=0}\Upsilon^k_i\,\frac{\partial
\bar{\Alpha}\vphantom{\Alpha}^p_{j\,q}}{\partial x^k}-\sum^3_{k=0}
\Upsilon^k_j\,\frac{\partial\bar{\Alpha}\vphantom{\Alpha}^p_{i\,q}}
{\partial x^k}+\sum^2_{h=1}\left(\bar{\Alpha}\vphantom{\Alpha}^p_{i\,h}
\,\bar{\Alpha}\vphantom{\Alpha}^{\!h}_{j\,q}-\bar{\Alpha}
\vphantom{\Alpha}^p_{j\,h}\,\bar{\Alpha}\vphantom{\Alpha}^{\!h}_{i\,q}
\right)\,-\\
-\sum^J_{P=1}\dsize\msum{2}\Sb i_1,\,\ldots,\,i_\alpha\\
j_1,\,\ldots,\,j_\beta\\ \bar i_1,\,\ldots,\,\bar i_\nu\\ 
\bar j_1,\,\ldots,\,\bar j_\gamma\endSb
\msum{3}\Sb h_1,\,\ldots,\,h_m\\ k_1,\,\ldots,\,k_n\endSb
\left(\,\shave{\sum^\alpha_{\mu=1}\sum^2_{v_\mu=1}}\Alpha^{i_\mu}_{i\,v_\mu}
\ S^{\,i_1\ldots\,v_\mu\,\ldots\,i_\alpha\,\bar i_1\ldots\,\bar i_\nu
\,h_1\ldots\,h_m}_{j_1\ldots\,\ldots\,\ldots\,j_\beta\,\bar j_1\ldots\,
\bar j_\gamma\,k_1\ldots\,k_n}[P]\,-\right.\\
\vspace{0.5ex plus 0.5ex minus 0.5ex}
-\sum^\beta_{\mu=1}\sum^2_{w_\mu=1}\Alpha^{w_\mu}_{i\,j_\mu}
\ S^{\,i_1\ldots\,\ldots\,\ldots\,i_\alpha\,\bar i_1\ldots\,\bar i_\nu
\,h_1\ldots\,h_m}_{j_1\ldots\,w_\mu\,\ldots\,j_\beta\,\bar j_1\ldots\,
\bar j_\gamma\,k_1\ldots\,k_n}[P]
+\sum^\nu_{\mu=1}\sum^2_{v_\mu=1}
\bar{\Alpha}\vphantom{\Alpha}^{\bar i_\mu}_{i\,v_\mu}\,\times\\
\vspace{0.5ex plus 0.5ex minus 0.5ex}
\times\,S^{\,i_1\ldots\,i_\alpha\,\bar i_1\ldots\,v_\mu\,\ldots\,
\bar i_\nu\,h_1\ldots\,h_m}_{j_1\ldots\,j_\beta\,\bar j_1\ldots\,\ldots
\,\ldots\,\bar j_\gamma\,k_1\ldots\,k_n}[P]-
\sum^\gamma_{\mu=1}\sum^2_{w_\mu=1}
\bar{\Alpha}\vphantom{\Alpha}^{w_\mu}_{i\,\bar j_\mu}\,\times\\
\vspace{0.5ex plus 0.5ex minus 0.5ex}
\times\,S^{i_1\ldots\,i_\alpha\,\bar i_1\ldots\,\ldots\,\ldots\,\bar i_\nu
\,h_1\ldots\,h_m}_{j_1\ldots\,j_\beta\,\bar j_1\ldots\,w_\mu\,\ldots\,
\bar j_\gamma\,k_1\ldots\,k_n}[P]
+\sum^m_{\mu=1}\sum^3_{v_\mu=0}\Gamma^{h_\mu}_{i\,v_\mu}\,\times\kern 4em\\
\vspace{0.5ex plus 0.5ex minus 0.5ex}
\times\,S^{i_1\ldots\,i_\alpha\,\bar i_1\ldots\,\bar i_\nu
\,h_1\ldots\,v_\mu\,\ldots\,h_m}_{j_1\ldots\,j_\beta\,\bar j_1\ldots\,
\bar j_\gamma\,k_1\ldots\,\ldots\,\ldots\,k_n}[P]-
\sum^n_{\mu=1}\sum^3_{w_\mu=0}\Gamma^{w_\mu}_{i\,k_\mu}\,\times\\
\endgathered\qquad
\mytag{6.26}\\
\gathered
\left.\vphantom{\shave{\sum^\alpha_{\mu=1}\sum^2_{v_\mu=1}}}
\times\,S^{i_1\ldots\,i_\alpha\,\bar i_1\ldots\,\bar i_\nu\,
h_1\ldots\,\ldots\,\ldots\,h_m}_{j_1\ldots\,j_\beta\,\bar j_1\ldots\,
\bar j_\gamma\,k_1\ldots\,w_\mu\,\ldots\,k_n}[P]\right)
\frac{\partial\bar{\Alpha}\vphantom{\Alpha}^p_{j\,q}}{\partial
S^{\,i_1\ldots\,i_\alpha\,\bar i_1\ldots\,\bar i_\nu\,h_1\ldots\,
h_m}_{j_1\ldots\,j_\beta\,\bar j_1\ldots\,\bar j_\gamma\,k_1\ldots
\,k_n}[P]}\,-\\
\vspace{0.132ex plus 0.5ex minus 0.7ex}
-\sum^J_{P=1}\dsize\msum{2}\Sb i_1,\,\ldots,\,i_\alpha\\
j_1,\,\ldots,\,j_\beta\\ \bar i_1,\,\ldots,\,\bar i_\nu\\ 
\bar j_1,\,\ldots,\,\bar j_\gamma\endSb
\msum{3}\Sb h_1,\,\ldots,\,h_m\\ k_1,\,\ldots,\,k_n\endSb
\left(\,\shave{\sum^\nu_{\mu=1}\sum^2_{v_\mu=1}}\Alpha^{i_\mu}_{i\,v_\mu}
\ \overline{S^{\,\bar i_1\ldots\,\bar i_\alpha\,i_1\ldots\,v_\mu\,\ldots\,
i_\nu h_1\ldots\,h_m}_{\bar j_1\ldots\,\bar j_\beta\,j_1\ldots\,\ldots\,
\ldots\,j_\gamma k_1\ldots\,k_n}[P]}\,-\right.\\
\vspace{0.132ex plus 0.5ex minus 0.7ex}
-\sum^\gamma_{\mu=1}\sum^2_{w_\mu=1}\Alpha^{w_\mu}_{i\,j_\mu}
\ \overline{S^{\,\bar i_1\ldots\,\bar i_\alpha\,i_1\ldots\,\ldots\,\ldots
\,i_\nu\,h_1\ldots\,h_m}_{\bar j_1\ldots\,\bar j_\beta\,j_1\ldots\,w_\mu
\,\ldots\,j_\gamma\,k_1\ldots\,k_n}[P]}+\sum^\alpha_{\mu=1}\sum^2_{v_\mu=1}
\bar{\Alpha}\vphantom{\Alpha}^{\bar i_\mu}_{i\,v_\mu}\,\times\\
\vspace{0.132ex plus 0.5ex minus 0.7ex}
\times\,\overline{S^{\,\bar i_1\ldots\,v_\mu\,\ldots\,\bar i_\alpha\,i_1
\ldots\,i_\nu\,h_1\ldots\,h_m}_{\bar j_1\ldots\,\ldots\,\ldots\,
\bar j_\beta\,j_1\ldots\,j_\gamma\,k_1\ldots\,k_n}[P]}
-\sum^\beta_{\mu=1}\sum^2_{w_\mu=1}
\bar{\Alpha}\vphantom{\Alpha}^{w_\mu}_{i\,\bar j_\mu}\,\times\\
\vspace{0.132ex plus 0.5ex minus 0.7ex}
\times\,\overline{S^{\,\bar i_1\ldots\,\ldots\,\ldots\,\bar i_\alpha\,
i_1\ldots\,i_\nu\,h_1\ldots\,h_m}_{\bar j_1\ldots\,w_\mu\,\ldots\,
\bar j_\beta\,j_1\ldots\,j_\gamma\,k_1\ldots\,k_n}[P]}+
\sum^m_{\mu=1}\sum^3_{v_\mu=0}\Gamma^{h_\mu}_{i\,v_\mu}\,\times\kern 4em\\
\times\,\overline{S^{\bar i_1\ldots\,\bar i_\alpha\,i_1\ldots\,i_\nu\,
h_1\ldots\,v_\mu\,\ldots\,h_m}_{\bar j_1\ldots\,\bar j_\beta\,j_1
\ldots\,j_\gamma\,k_1\ldots\,\ldots\,\ldots\,k_n}[P]}-
\sum^n_{\mu=1}\sum^3_{w_\mu=0}\Gamma^{w_\mu}_{i\,k_\mu}\,\times\\
\vspace{0.132ex plus 0.5ex minus 0.7ex}
\left.\vphantom{\shave{\sum^\nu_{\mu=1}\sum^2_{v_\mu=1}}}
\times\,\overline{S^{\bar i_1\ldots\,\bar i_\alpha\,i_1\ldots\,i_\nu
\,h_1\ldots\,\ldots\,\ldots\,h_m}_{\bar j_1\ldots\,\bar j_\beta
j_1\ldots\,j_\gamma\,k_1\ldots\,w_\mu\,\ldots\,k_n}[P]}
\right)\frac{\partial\bar{\Alpha}\vphantom{\Alpha}^p_{j\,q}}
{\partial\overline{S^{\,i_1\ldots\,i_\alpha\,\bar i_1\ldots\,\bar i_\nu
\,h_1\ldots\,h_m}_{j_1\ldots\,j_\beta\,\bar j_1\ldots\,\bar j_\gamma\,
k_1\ldots\,k_n}[P]}\,}\,+\\
\vspace{0.132ex plus 0.5ex minus 0.7ex}
+\sum^J_{P=1}\dsize\msum{2}\Sb i_1,\,\ldots,\,i_\alpha\\
j_1,\,\ldots,\,j_\beta\\ \bar i_1,\,\ldots,\,\bar i_\nu\\ 
\bar j_1,\,\ldots,\,\bar j_\gamma\endSb
\msum{3}\Sb h_1,\,\ldots,\,h_m\\ k_1,\,\ldots,\,k_n\endSb
\left(\,\shave{\sum^\alpha_{\mu=1}\sum^2_{v_\mu=1}}\Alpha^{i_\mu}_{j\,v_\mu}
\ S^{\,i_1\ldots\,v_\mu\,\ldots\,i_\alpha\,\bar i_1\ldots\,\bar i_\nu
\,h_1\ldots\,h_m}_{j_1\ldots\,\ldots\,\ldots\,j_\beta\,\bar j_1\ldots\,
\bar j_\gamma\,k_1\ldots\,k_n}[P]\,-\right.\\
\vspace{0.132ex plus 0.5ex minus 0.7ex}
-\sum^\beta_{\mu=1}\sum^2_{w_\mu=1}\Alpha^{w_\mu}_{j\,j_\mu}
\ S^{\,i_1\ldots\,\ldots\,\ldots\,i_\alpha\,\bar i_1\ldots\,\bar i_\nu
\,h_1\ldots\,h_m}_{j_1\ldots\,w_\mu\,\ldots\,j_\beta\,\bar j_1\ldots\,
\bar j_\gamma\,k_1\ldots\,k_n}[P]
+\sum^\nu_{\mu=1}\sum^2_{v_\mu=1}
\bar{\Alpha}\vphantom{\Alpha}^{\bar i_\mu}_{j\,v_\mu}\,\times\\
\vspace{0.132ex plus 0.5ex minus 0.7ex}
\times\,S^{\,i_1\ldots\,i_\alpha\,\bar i_1\ldots\,v_\mu\,\ldots\,
\bar i_\nu\,h_1\ldots\,h_m}_{j_1\ldots\,j_\beta\,\bar j_1\ldots\,\ldots
\,\ldots\,\bar j_\gamma\,k_1\ldots\,k_n}[P]-
\sum^\gamma_{\mu=1}\sum^2_{w_\mu=1}
\bar{\Alpha}\vphantom{\Alpha}^{w_\mu}_{j\,\bar j_\mu}\,\times\\
\vspace{0.132ex plus 0.5ex minus 0.7ex}
\times\,S^{i_1\ldots\,i_\alpha\,\bar i_1\ldots\,\ldots\,\ldots\,\bar i_\nu
\,h_1\ldots\,h_m}_{j_1\ldots\,j_\beta\,\bar j_1\ldots\,w_\mu\,\ldots\,
\bar j_\gamma\,k_1\ldots\,k_n}[P]
+\sum^m_{\mu=1}\sum^3_{v_\mu=0}\Gamma^{h_\mu}_{j\,v_\mu}\,\times\kern 4em\\
\vspace{0.132ex plus 0.5ex minus 0.7ex}
\times\,S^{i_1\ldots\,i_\alpha\,\bar i_1\ldots\,\bar i_\nu
\,h_1\ldots\,v_\mu\,\ldots\,h_m}_{j_1\ldots\,j_\beta\,\bar j_1\ldots\,
\bar j_\gamma\,k_1\ldots\,\ldots\,\ldots\,k_n}[P]-
\sum^n_{\mu=1}\sum^3_{w_\mu=0}\Gamma^{w_\mu}_{j\,k_\mu}\,\times\\
\vspace{0.132ex plus 0.5ex minus 0.7ex}
\left.\vphantom{\shave{\sum^\alpha_{\mu=1}\sum^2_{v_\mu=1}}}
\times\,S^{i_1\ldots\,i_\alpha\,\bar i_1\ldots\,\bar i_\nu\,
h_1\ldots\,\ldots\,\ldots\,h_m}_{j_1\ldots\,j_\beta\,\bar j_1\ldots\,
\bar j_\gamma\,k_1\ldots\,w_\mu\,\ldots\,k_n}[P]\right)
\frac{\partial\bar{\Alpha}\vphantom{\Alpha}^p_{i\,q}}{\partial
S^{\,i_1\ldots\,i_\alpha\,\bar i_1\ldots\,\bar i_\nu\,h_1\ldots\,
h_m}_{j_1\ldots\,j_\beta\,\bar j_1\ldots\,\bar j_\gamma\,k_1\ldots
\,k_n}[P]}\,+\\
\vspace{0.132ex plus 0.5ex minus 0.7ex}
+\sum^J_{P=1}\dsize\msum{2}\Sb i_1,\,\ldots,\,i_\alpha\\
j_1,\,\ldots,\,j_\beta\\ \bar i_1,\,\ldots,\,\bar i_\nu\\ 
\bar j_1,\,\ldots,\,\bar j_\gamma\endSb
\msum{3}\Sb h_1,\,\ldots,\,h_m\\ k_1,\,\ldots,\,k_n\endSb
\left(\,\shave{\sum^\nu_{\mu=1}\sum^2_{v_\mu=1}}\Alpha^{i_\mu}_{j\,v_\mu}
\ \overline{S^{\,\bar i_1\ldots\,\bar i_\alpha\,i_1\ldots\,v_\mu\,\ldots\,
i_\nu h_1\ldots\,h_m}_{\bar j_1\ldots\,\bar j_\beta\,j_1\ldots\,\ldots\,
\ldots\,j_\gamma k_1\ldots\,k_n}[P]}\,-\right.
\endgathered\kern 4em\\
\gathered
%\vspace{0.5ex plus 0.5ex minus 0.5ex}
-\sum^\gamma_{\mu=1}\sum^2_{w_\mu=1}\Alpha^{w_\mu}_{j\,j_\mu}
\ \overline{S^{\,\bar i_1\ldots\,\bar i_\alpha\,i_1\ldots\,\ldots\,\ldots
\,i_\nu\,h_1\ldots\,h_m}_{\bar j_1\ldots\,\bar j_\beta\,j_1\ldots\,w_\mu
\,\ldots\,j_\gamma\,k_1\ldots\,k_n}[P]}+\sum^\alpha_{\mu=1}\sum^2_{v_\mu=1}
\bar{\Alpha}\vphantom{\Alpha}^{\bar i_\mu}_{j\,v_\mu}\,\times\\
\vspace{0.5ex plus 0.5ex minus 0.5ex}
\times\,\overline{S^{\,\bar i_1\ldots\,v_\mu\,\ldots\,\bar i_\alpha\,i_1
\ldots\,i_\nu\,h_1\ldots\,h_m}_{\bar j_1\ldots\,\ldots\,\ldots\,
\bar j_\beta\,j_1\ldots\,j_\gamma\,k_1\ldots\,k_n}[P]}
-\sum^\beta_{\mu=1}\sum^2_{w_\mu=1}
\bar{\Alpha}\vphantom{\Alpha}^{w_\mu}_{j\,\bar j_\mu}\,\times\\
\vspace{0.5ex plus 0.5ex minus 0.5ex}
\times\,\overline{S^{\,\bar i_1\ldots\,\ldots\,\ldots\,\bar i_\alpha\,
i_1\ldots\,i_\nu\,h_1\ldots\,h_m}_{\bar j_1\ldots\,w_\mu\,\ldots\,
\bar j_\beta\,j_1\ldots\,j_\gamma\,k_1\ldots\,k_n}[P]}+
\sum^m_{\mu=1}\sum^3_{v_\mu=0}\Gamma^{h_\mu}_{j\,v_\mu}\,\times\kern 4em\\
\vspace{0.5ex plus 0.5ex minus 0.5ex}
\times\,\overline{S^{\bar i_1\ldots\,\bar i_\alpha\,i_1\ldots\,i_\nu\,
h_1\ldots\,v_\mu\,\ldots\,h_m}_{\bar j_1\ldots\,\bar j_\beta\,j_1
\ldots\,j_\gamma\,k_1\ldots\,\ldots\,\ldots\,k_n}[P]}-
\sum^n_{\mu=1}\sum^3_{w_\mu=0}\Gamma^{w_\mu}_{j\,k_\mu}\,\times\\
\vspace{0.5ex plus 0.5ex minus 0.5ex}
\left.\vphantom{\shave{\sum^\nu_{\mu=1}\sum^2_{v_\mu=1}}}
\times\,\overline{S^{\bar i_1\ldots\,\bar i_\alpha\,i_1\ldots\,i_\nu
\,h_1\ldots\,\ldots\,\ldots\,h_m}_{\bar j_1\ldots\,\bar j_\beta
j_1\ldots\,j_\gamma\,k_1\ldots\,w_\mu\,\ldots\,k_n}[P]}
\right)\frac{\partial\bar{\Alpha}\vphantom{\Alpha}^p_{i\,q}}
{\partial\overline{S^{\,i_1\ldots\,i_\alpha\,\bar i_1\ldots\,\bar i_\nu
\,h_1\ldots\,h_m}_{j_1\ldots\,j_\beta\,\bar j_1\ldots\,\bar j_\gamma\,
k_1\ldots\,k_n}[P]}\,}-\sum^3_{k=0}c^{\,k}_{ij}\,\bar{\Alpha}
\vphantom{\Alpha}^p_{kq}.
\endgathered\kern 4em
\endgather
$$
The components of the third curvature spin-tensor $\bold R$ are given
by the formula
$$
\allowdisplaybreaks
\gather
\gathered
R^p_{qij}=\sum^3_{k=0}\Upsilon^k_i\,\frac{\partial
\Gamma^p_{\!j\,q}}{\partial x^k}-\sum^3_{k=0}\Upsilon^k_j
\,\frac{\partial\Gamma^p_{\!i\,q}}{\partial x^k}
+\sum^2_{h=1}\left(\Gamma^p_{\!i\,h}\,\Gamma^h_{\!j\,q}
-\Gamma^p_{\!j\,h}\,\Gamma^h_{\!i\,q}\right)\,-\\
-\sum^J_{P=1}\dsize\msum{2}\Sb i_1,\,\ldots,\,i_\alpha\\
j_1,\,\ldots,\,j_\beta\\ \bar i_1,\,\ldots,\,\bar i_\nu\\ 
\bar j_1,\,\ldots,\,\bar j_\gamma\endSb
\msum{3}\Sb h_1,\,\ldots,\,h_m\\ k_1,\,\ldots,\,k_n\endSb
\left(\,\shave{\sum^\alpha_{\mu=1}\sum^2_{v_\mu=1}}\Alpha^{i_\mu}_{i\,v_\mu}
\ S^{\,i_1\ldots\,v_\mu\,\ldots\,i_\alpha\,\bar i_1\ldots\,\bar i_\nu
\,h_1\ldots\,h_m}_{j_1\ldots\,\ldots\,\ldots\,j_\beta\,\bar j_1\ldots\,
\bar j_\gamma\,k_1\ldots\,k_n}[P]\,-\right.\\
\vspace{0.5ex plus 0.5ex minus 0.5ex}
-\sum^\beta_{\mu=1}\sum^2_{w_\mu=1}\Alpha^{w_\mu}_{i\,j_\mu}
\ S^{\,i_1\ldots\,\ldots\,\ldots\,i_\alpha\,\bar i_1\ldots\,\bar i_\nu
\,h_1\ldots\,h_m}_{j_1\ldots\,w_\mu\,\ldots\,j_\beta\,\bar j_1\ldots\,
\bar j_\gamma\,k_1\ldots\,k_n}[P]
+\sum^\nu_{\mu=1}\sum^2_{v_\mu=1}
\bar{\Alpha}\vphantom{\Alpha}^{\bar i_\mu}_{i\,v_\mu}\,\times\\
\vspace{0.5ex plus 0.5ex minus 0.5ex}
\times\,S^{\,i_1\ldots\,i_\alpha\,\bar i_1\ldots\,v_\mu\,\ldots\,
\bar i_\nu\,h_1\ldots\,h_m}_{j_1\ldots\,j_\beta\,\bar j_1\ldots\,\ldots
\,\ldots\,\bar j_\gamma\,k_1\ldots\,k_n}[P]-
\sum^\gamma_{\mu=1}\sum^2_{w_\mu=1}
\bar{\Alpha}\vphantom{\Alpha}^{w_\mu}_{i\,\bar j_\mu}\,\times\\
\vspace{0.5ex plus 0.5ex minus 0.5ex}
\times\,S^{i_1\ldots\,i_\alpha\,\bar i_1\ldots\,\ldots\,\ldots\,\bar i_\nu
\,h_1\ldots\,h_m}_{j_1\ldots\,j_\beta\,\bar j_1\ldots\,w_\mu\,\ldots\,
\bar j_\gamma\,k_1\ldots\,k_n}[P]
+\sum^m_{\mu=1}\sum^3_{v_\mu=0}\Gamma^{h_\mu}_{i\,v_\mu}\,\times\kern 4em\\
\vspace{0.5ex plus 0.5ex minus 0.5ex}
\times\,S^{i_1\ldots\,i_\alpha\,\bar i_1\ldots\,\bar i_\nu
\,h_1\ldots\,v_\mu\,\ldots\,h_m}_{j_1\ldots\,j_\beta\,\bar j_1\ldots\,
\bar j_\gamma\,k_1\ldots\,\ldots\,\ldots\,k_n}[P]-
\sum^n_{\mu=1}\sum^3_{w_\mu=0}\Gamma^{w_\mu}_{i\,k_\mu}\,\times\\
\vspace{0.5ex plus 0.5ex minus 0.5ex}
\left.\vphantom{\shave{\sum^\alpha_{\mu=1}\sum^2_{v_\mu=1}}}
\times\,S^{i_1\ldots\,i_\alpha\,\bar i_1\ldots\,\bar i_\nu\,
h_1\ldots\,\ldots\,\ldots\,h_m}_{j_1\ldots\,j_\beta\,\bar j_1\ldots\,
\bar j_\gamma\,k_1\ldots\,w_\mu\,\ldots\,k_n}[P]\right)
\frac{\partial\Gamma^p_{j\,q}}{\partial
S^{\,i_1\ldots\,i_\alpha\,\bar i_1\ldots\,\bar i_\nu\,h_1\ldots\,
h_m}_{j_1\ldots\,j_\beta\,\bar j_1\ldots\,\bar j_\gamma\,k_1\ldots
\,k_n}[P]}\,-\\
\vspace{0.5ex plus 0.5ex minus 0.5ex}
-\sum^J_{P=1}\dsize\msum{2}\Sb i_1,\,\ldots,\,i_\alpha\\
j_1,\,\ldots,\,j_\beta\\ \bar i_1,\,\ldots,\,\bar i_\nu\\ 
\bar j_1,\,\ldots,\,\bar j_\gamma\endSb
\msum{3}\Sb h_1,\,\ldots,\,h_m\\ k_1,\,\ldots,\,k_n\endSb
\left(\,\shave{\sum^\nu_{\mu=1}\sum^2_{v_\mu=1}}\Alpha^{i_\mu}_{i\,v_\mu}
\ \overline{S^{\,\bar i_1\ldots\,\bar i_\alpha\,i_1\ldots\,v_\mu\,\ldots\,
i_\nu h_1\ldots\,h_m}_{\bar j_1\ldots\,\bar j_\beta\,j_1\ldots\,\ldots\,
\ldots\,j_\gamma k_1\ldots\,k_n}[P]}\,-\right.
\endgathered\qquad
\mytag{6.27}\\
\gathered
-\sum^\gamma_{\mu=1}\sum^2_{w_\mu=1}\Alpha^{w_\mu}_{i\,j_\mu}
\ \overline{S^{\,\bar i_1\ldots\,\bar i_\alpha\,i_1\ldots\,\ldots\,\ldots
\,i_\nu\,h_1\ldots\,h_m}_{\bar j_1\ldots\,\bar j_\beta\,j_1\ldots\,w_\mu
\,\ldots\,j_\gamma\,k_1\ldots\,k_n}[P]}+\sum^\alpha_{\mu=1}\sum^2_{v_\mu=1}
\bar{\Alpha}\vphantom{\Alpha}^{\bar i_\mu}_{i\,v_\mu}\,\times\\
\vspace{0.575ex plus 0.5ex minus 0.5ex}
\times\,\overline{S^{\,\bar i_1\ldots\,v_\mu\,\ldots\,\bar i_\alpha\,i_1
\ldots\,i_\nu\,h_1\ldots\,h_m}_{\bar j_1\ldots\,\ldots\,\ldots\,
\bar j_\beta\,j_1\ldots\,j_\gamma\,k_1\ldots\,k_n}[P]}
-\sum^\beta_{\mu=1}\sum^2_{w_\mu=1}
\bar{\Alpha}\vphantom{\Alpha}^{w_\mu}_{i\,\bar j_\mu}\,\times\\
\vspace{0.575ex plus 0.5ex minus 0.5ex}
\times\,\overline{S^{\,\bar i_1\ldots\,\ldots\,\ldots\,\bar i_\alpha\,
i_1\ldots\,i_\nu\,h_1\ldots\,h_m}_{\bar j_1\ldots\,w_\mu\,\ldots\,
\bar j_\beta\,j_1\ldots\,j_\gamma\,k_1\ldots\,k_n}[P]}+
\sum^m_{\mu=1}\sum^3_{v_\mu=0}\Gamma^{h_\mu}_{i\,v_\mu}\,\times\kern 4em\\
\times\,\overline{S^{\bar i_1\ldots\,\bar i_\alpha\,i_1\ldots\,i_\nu\,
h_1\ldots\,v_\mu\,\ldots\,h_m}_{\bar j_1\ldots\,\bar j_\beta\,j_1
\ldots\,j_\gamma\,k_1\ldots\,\ldots\,\ldots\,k_n}[P]}-
\sum^n_{\mu=1}\sum^3_{w_\mu=0}\Gamma^{w_\mu}_{i\,k_\mu}\,\times\\
\vspace{0.575ex plus 0.5ex minus 0.5ex}
\left.\vphantom{\shave{\sum^\nu_{\mu=1}\sum^2_{v_\mu=1}}}
\times\,\overline{S^{\bar i_1\ldots\,\bar i_\alpha\,i_1\ldots\,i_\nu
\,h_1\ldots\,\ldots\,\ldots\,h_m}_{\bar j_1\ldots\,\bar j_\beta
j_1\ldots\,j_\gamma\,k_1\ldots\,w_\mu\,\ldots\,k_n}[P]}
\right)\frac{\partial\Gamma^p_{j\,q}}
{\partial\overline{S^{\,i_1\ldots\,i_\alpha\,\bar i_1\ldots\,\bar i_\nu
\,h_1\ldots\,h_m}_{j_1\ldots\,j_\beta\,\bar j_1\ldots\,\bar j_\gamma\,
k_1\ldots\,k_n}[P]}\,}\,+\\
\vspace{0.575ex plus 0.5ex minus 0.5ex}
+\sum^J_{P=1}\dsize\msum{2}\Sb i_1,\,\ldots,\,i_\alpha\\
j_1,\,\ldots,\,j_\beta\\ \bar i_1,\,\ldots,\,\bar i_\nu\\ 
\bar j_1,\,\ldots,\,\bar j_\gamma\endSb
\msum{3}\Sb h_1,\,\ldots,\,h_m\\ k_1,\,\ldots,\,k_n\endSb
\left(\,\shave{\sum^\alpha_{\mu=1}\sum^2_{v_\mu=1}}\Alpha^{i_\mu}_{j\,v_\mu}
\ S^{\,i_1\ldots\,v_\mu\,\ldots\,i_\alpha\,\bar i_1\ldots\,\bar i_\nu
\,h_1\ldots\,h_m}_{j_1\ldots\,\ldots\,\ldots\,j_\beta\,\bar j_1\ldots\,
\bar j_\gamma\,k_1\ldots\,k_n}[P]\,-\right.\\
\vspace{0.575ex plus 0.5ex minus 0.5ex}
-\sum^\beta_{\mu=1}\sum^2_{w_\mu=1}\Alpha^{w_\mu}_{j\,j_\mu}
\ S^{\,i_1\ldots\,\ldots\,\ldots\,i_\alpha\,\bar i_1\ldots\,\bar i_\nu
\,h_1\ldots\,h_m}_{j_1\ldots\,w_\mu\,\ldots\,j_\beta\,\bar j_1\ldots\,
\bar j_\gamma\,k_1\ldots\,k_n}[P]
+\sum^\nu_{\mu=1}\sum^2_{v_\mu=1}
\bar{\Alpha}\vphantom{\Alpha}^{\bar i_\mu}_{j\,v_\mu}\,\times\\
\vspace{0.575ex plus 0.5ex minus 0.5ex}
\times\,S^{\,i_1\ldots\,i_\alpha\,\bar i_1\ldots\,v_\mu\,\ldots\,
\bar i_\nu\,h_1\ldots\,h_m}_{j_1\ldots\,j_\beta\,\bar j_1\ldots\,\ldots
\,\ldots\,\bar j_\gamma\,k_1\ldots\,k_n}[P]-
\sum^\gamma_{\mu=1}\sum^2_{w_\mu=1}
\bar{\Alpha}\vphantom{\Alpha}^{w_\mu}_{j\,\bar j_\mu}\,\times\\
\vspace{0.575ex plus 0.5ex minus 0.5ex}
\times\,S^{i_1\ldots\,i_\alpha\,\bar i_1\ldots\,\ldots\,\ldots\,\bar i_\nu
\,h_1\ldots\,h_m}_{j_1\ldots\,j_\beta\,\bar j_1\ldots\,w_\mu\,\ldots\,
\bar j_\gamma\,k_1\ldots\,k_n}[P]
+\sum^m_{\mu=1}\sum^3_{v_\mu=0}\Gamma^{h_\mu}_{j\,v_\mu}\,\times\kern 4em\\
\vspace{0.575ex plus 0.5ex minus 0.5ex}
\times\,S^{i_1\ldots\,i_\alpha\,\bar i_1\ldots\,\bar i_\nu
\,h_1\ldots\,v_\mu\,\ldots\,h_m}_{j_1\ldots\,j_\beta\,\bar j_1\ldots\,
\bar j_\gamma\,k_1\ldots\,\ldots\,\ldots\,k_n}[P]-
\sum^n_{\mu=1}\sum^3_{w_\mu=0}\Gamma^{w_\mu}_{j\,k_\mu}\,\times\\
\vspace{0.575ex plus 0.5ex minus 0.5ex}
\left.\vphantom{\shave{\sum^\alpha_{\mu=1}\sum^2_{v_\mu=1}}}
\times\,S^{i_1\ldots\,i_\alpha\,\bar i_1\ldots\,\bar i_\nu\,
h_1\ldots\,\ldots\,\ldots\,h_m}_{j_1\ldots\,j_\beta\,\bar j_1\ldots\,
\bar j_\gamma\,k_1\ldots\,w_\mu\,\ldots\,k_n}[P]\right)
\frac{\partial\Gamma^p_{i\,q}}{\partial
S^{\,i_1\ldots\,i_\alpha\,\bar i_1\ldots\,\bar i_\nu\,h_1\ldots\,
h_m}_{j_1\ldots\,j_\beta\,\bar j_1\ldots\,\bar j_\gamma\,k_1\ldots
\,k_n}[P]}\,+\\
\vspace{0.575ex plus 0.5ex minus 0.5ex}
+\sum^J_{P=1}\dsize\msum{2}\Sb i_1,\,\ldots,\,i_\alpha\\
j_1,\,\ldots,\,j_\beta\\ \bar i_1,\,\ldots,\,\bar i_\nu\\ 
\bar j_1,\,\ldots,\,\bar j_\gamma\endSb
\msum{3}\Sb h_1,\,\ldots,\,h_m\\ k_1,\,\ldots,\,k_n\endSb
\left(\,\shave{\sum^\nu_{\mu=1}\sum^2_{v_\mu=1}}\Alpha^{i_\mu}_{j\,v_\mu}
\ \overline{S^{\,\bar i_1\ldots\,\bar i_\alpha\,i_1\ldots\,v_\mu\,\ldots\,
i_\nu h_1\ldots\,h_m}_{\bar j_1\ldots\,\bar j_\beta\,j_1\ldots\,\ldots\,
\ldots\,j_\gamma k_1\ldots\,k_n}[P]}\,-\right.\\
\vspace{0.575ex plus 0.5ex minus 0.5ex}
-\sum^\gamma_{\mu=1}\sum^2_{w_\mu=1}\Alpha^{w_\mu}_{j\,j_\mu}
\ \overline{S^{\,\bar i_1\ldots\,\bar i_\alpha\,i_1\ldots\,\ldots\,\ldots
\,i_\nu\,h_1\ldots\,h_m}_{\bar j_1\ldots\,\bar j_\beta\,j_1\ldots\,w_\mu
\,\ldots\,j_\gamma\,k_1\ldots\,k_n}[P]}+\sum^\alpha_{\mu=1}\sum^2_{v_\mu=1}
\bar{\Alpha}\vphantom{\Alpha}^{\bar i_\mu}_{j\,v_\mu}\,\times\\
\vspace{0.575ex plus 0.5ex minus 0.5ex}
\times\,\overline{S^{\,\bar i_1\ldots\,v_\mu\,\ldots\,\bar i_\alpha\,i_1
\ldots\,i_\nu\,h_1\ldots\,h_m}_{\bar j_1\ldots\,\ldots\,\ldots\,
\bar j_\beta\,j_1\ldots\,j_\gamma\,k_1\ldots\,k_n}[P]}
-\sum^\beta_{\mu=1}\sum^2_{w_\mu=1}
\bar{\Alpha}\vphantom{\Alpha}^{w_\mu}_{j\,\bar j_\mu}\,\times\\
\endgathered\kern 4em\\
\gathered
\times\,\overline{S^{\,\bar i_1\ldots\,\ldots\,\ldots\,\bar i_\alpha\,
i_1\ldots\,i_\nu\,h_1\ldots\,h_m}_{\bar j_1\ldots\,w_\mu\,\ldots\,
\bar j_\beta\,j_1\ldots\,j_\gamma\,k_1\ldots\,k_n}[P]}+
\sum^m_{\mu=1}\sum^3_{v_\mu=0}\Gamma^{h_\mu}_{j\,v_\mu}\,\times\kern 4em\\
\vspace{0.5ex plus 0.5ex minus 0.5ex}
\times\,\overline{S^{\bar i_1\ldots\,\bar i_\alpha\,i_1\ldots\,i_\nu\,
h_1\ldots\,v_\mu\,\ldots\,h_m}_{\bar j_1\ldots\,\bar j_\beta\,j_1
\ldots\,j_\gamma\,k_1\ldots\,\ldots\,\ldots\,k_n}[P]}-
\sum^n_{\mu=1}\sum^3_{w_\mu=0}\Gamma^{w_\mu}_{j\,k_\mu}\,\times\\
\vspace{0.5ex plus 0.5ex minus 0.5ex}
\left.\vphantom{\shave{\sum^\nu_{\mu=1}\sum^2_{v_\mu=1}}}
\times\,\overline{S^{\bar i_1\ldots\,\bar i_\alpha\,i_1\ldots\,i_\nu
\,h_1\ldots\,\ldots\,\ldots\,h_m}_{\bar j_1\ldots\,\bar j_\beta
j_1\ldots\,j_\gamma\,k_1\ldots\,w_\mu\,\ldots\,k_n}[P]}
\right)\frac{\partial\Gamma^p_{i\,q}}
{\partial\overline{S^{\,i_1\ldots\,i_\alpha\,\bar i_1\ldots\,\bar i_\nu
\,h_1\ldots\,h_m}_{j_1\ldots\,j_\beta\,\bar j_1\ldots\,\bar j_\gamma\,
k_1\ldots\,k_n}[P]}\,}-\sum^3_{k=0}c^{\,k}_{ij}\,\Gamma^p_{kq}.
\endgathered\kern 4em
\endgather
$$
Like in \mythetag{6.25}, in the above formulas \mythetag{6.26} and
\mythetag{6.27} for each $P$ we implicitly assume that $\alpha=\alpha_P$, 
$\beta=\beta_P$, $\nu=\nu_P$, $\gamma=\gamma_P$, $m=m_P$, $n=n_P$.
Moreover, like in \mythetag{6.22} the quantities $c^{\,k}_{ij}$ in the 
last terms of \mythetag{6.25}, \mythetag{6.26}, and \mythetag{6.27} are
taken from \mythetag{1.14} or from \mythetag{1.15}.
\head
7. Acknowledgments.
\endhead
    I am grateful to E.~G.~Neufeld who gave me the book \mycite{3}. I am
also grateful to him for discussing the concept of spinors when I bagan
reading this book.

\enddocument
\end